\documentclass[a4paper,11pt,oneside]{article}

\usepackage{amsmath,amsthm,amsfonts,amssymb,amscd}
\usepackage{graphicx}
\usepackage{mathtools}
\usepackage{pgf,tikz}
\usepackage{mathrsfs}
\usepackage{float}
\usepackage{csvsimple}
\usepackage{enumerate}
\usepackage{booktabs}
\usepackage{multirow}
\usepackage{siunitx}
\usepackage{pgfplots}
\usepackage[noend]{algpseudocode}
\usepackage{graphicx,bm,xcolor}
\usepackage{algorithm}
\usepackage{pdfpages}
\usepackage{stmaryrd}
\usetikzlibrary{arrows}
\usepackage{caption}
\usepackage[T1]{fontenc}
\usepackage{t1enc}
\usepackage[polish,german,english]{babel}
\usepackage{verbatim} 
\usepackage{xcolor}

\usepackage{url}

\newif\ifMAKEPICS
\MAKEPICSfalse

\ifMAKEPICS
\usepackage[cleanup,subfolder]{gnuplottex}
\usepackage{xparse}

\ExplSyntaxOn
\DeclareExpandableDocumentCommand{\convertlen}{ O{cm} m }
{
	\dim_to_decimal_in_unit:nn { #2 } { 1 #1 } cm
}
\ExplSyntaxOff
\fi

\pgfplotsset{compat=1.16}

\usepackage{authblk}

\usepackage[most]{tcolorbox}

\usepackage[top=2cm, bottom=2cm, left=2cm, right=2cm]{geometry}

\theoremstyle{plain} 
\newtheorem{theorem}{Theorem}[section]

\newtheorem{corollary}[theorem]{Corollary}
\newtheorem{remark}[theorem]{Remark}

\theoremstyle{definition} %
\newtheorem{assumption}{Assumption}

\theoremstyle{remark} %

\newcommand{\ignore}[1]{}

\begin{document}
\title{Goal-Oriented Adaptive Space-Time Finite Element Methods\\for Regularized Parabolic p-Laplace Problems}
	\author[1,2]{B. Endtmayer}
	\author[3,4]{U. Langer}
	\author[3,4]{A. Schafelner}

\affil[1]{Leibniz Universit\"at Hannover,
	Institut f\"ur Angewandte Mathematik,
	AG Numerische Mathematik,
	Welfengarten 1, 30167 Hannover, Germany}
\affil[2]{Cluster of Excellence PhoenixD,
            (Photonics, Optics, and Engineering - Innovation Across Disciplines), 
        	Leibniz Universit\"at Hannover, Germany}
\affil[3]{Institute of Numerical Mathematics, Johannes Kepler University Linz, Altenbergerstr. 69, A-4040 Linz, Austria}
\affil[4]{Johann Radon Institute for Computational and Applied Mathematics, Altenbergerstr. 69, A-4040 Linz, Austria}

	\date{}
	
	\maketitle

\begin{abstract}
We consider goal-oriented adaptive space-time finite-element discretizations of the regularized 
parabolic p-Laplace problem on completely unstructured simplicial space-time meshes.
The adaptivity is driven by the dual-weighted residual (DWR) method since we are interested 
in an accurate computation of some possibly nonlinear functionals at the solution.
Such functionals represent goals in which engineers are often more interested than the solution itself.
The DWR method requires the numerical solution of a linear adjoint problem 
that provides the sensitivities for the mesh refinement. This can be done 
by means of the same full space-time finite element discretization as used 
for the primal non-linear problems.
The numerical experiments presented demonstrate that this goal-oriented, full space-time  
finite element solver efficiently provides accurate numerical results 
for different functionals.
\end{abstract}


\begin{msc}
35K61, 
65M50,  
65M60 
\end{msc}

\section{Introduction}
\label{Section:Introduction}

In this paper, we investigate goal-oriented adaptive space-time (GOAST) finite element discretizations of 
initial-boundary value problems (IBVP) for the  scalar regularized parabolic $p$-Laplace equation
\begin{equation}
\label{eqn:Parabolic_p-Laplace_CF} 
\partial_{t}u - \mathrm{div}_x((|\nabla_{x}u|^2 + \varepsilon^2)^{\frac{p-2}{2}}\nabla_{x}u) = f\;\text{in}\;Q,
\quad u=u_D:=0 \;\text{on}\; \Sigma, \quad u=u_0:=0 \;\text{on}\; \overline{\Sigma}_0
\end{equation}
on completely unstructured simplicial space-time meshes, where
$Q=\Omega\times(0,T) \subset \mathbb{R}^{d+1}$ denotes the space-time cylinder, 
$\Sigma=\partial{\Omega}\times(0,T)$ its lateral surface where the Dirichlet boundary condition $u_D$ is prescribed, 
$\Sigma_0 = \Omega\times\{0\}$ its bottom where the initial condition $u_0$ is given,
$\Omega \subset \mathbb{R}^{d}$ with the boundary $\partial{\Omega}$ denotes the spatial domain 
that is supposed to be bounded and Lipschitz, $d \in \{1,2,3\}$ is the spatial dimension, 
$T > 0$ is the final time,
the right-hand side $f$ is a given source driving the evolution process, the given power 
$p \in (1, \infty)$ characterizes the nonlinearity of the evolution, 
and $\varepsilon > 0$ is a positive regularization parameter.

In the limit case $\varepsilon = 0$, we speak about the original parabolic $p$-Laplace problem
that corresponds to the simple (linear) transient heat problem for $p=2$,
but that can degenerate otherwise.
As usual, $\nabla_x$,  $\text{div}_x$, and $\partial_t$ stand for the spatial gradient, 
the spatial divergence, and the partial time-derivative, respectively,
whereas 
$|\nabla_{x}u|=\big( \sum_{i=1}^d |\partial_i u|^2\big)^{1/2}$ denotes the Euclidean norm of 
the spatial gradient $\nabla_{x}u = (\partial_1 u,\ldots, \partial_d u)^\top{}$ of $u$.
We note that the regularized parabolic $p$-Laplace equation is sometimes 
written in the form $\partial_{t}u - \mathrm{div}_x((|\nabla_{x}u| + \varepsilon)^{p-2}\nabla_{x}u) = f$ instead of \eqref{eqn:Parabolic_p-Laplace_CF}; see, e.g., 
the book \cite{Roubicek:2013a}, in particular, Exercises~4.26 therein, or 
the recent
	 works \cite{Toulopoulos:2022a,kaltenbach2022error,HirnWollner:2020}.
There are many publications on 
the solvability analysis including investigations of the regularity of the solution 
\cite{Zeidler:1990a,DiBenedetto:1993a,Roubicek:2013a,CianchiMazya:2020a},
the numerical solution including  
discretization error estimates 
\cite{BarrettLui:1994SINUM.pdf}
and solvers for the nonlinear finite element equations 
\cite{LeePark:2022a},
and applications of elliptic (stationary) and parabolic (evolutionary) 
$p$-Laplace models in different disciplines
\cite{Diaz:1985a},
and the references therein.

The standard numerical 
technique for
the discretization of parabolic 
$p$-Laplace problems like \eqref{eqn:Parabolic_p-Laplace_CF} 
uses some time-stepping 
for the temporal discretization
in combination with a spatial finite element Galerkin discretization;
see, e.g., \cite{BarrettLui:1994SINUM.pdf}.
In contrast to this traditional approach, 
we here propose to use a completely unstructured space-time finite 
element method on simplicial space-time meshes 
that finally leads to
the
solution of one nonlinear system of 
finite element equations
instead of many smaller nonlinear systems 
in the case of implicit time-stepping. 
This space-time approach has successfully been used for linear parabolic
initial boundary value problems; 
see, e.g.\ the survey article \cite{SteinbachYang:2019a}.
Concerning the parabolic $p$-Laplace initial-boundary value problem,
we are only aware of the 
very recent publication
\cite{Toulopoulos:2022a}
by Toulopoulos who 
derived consistent space-time finite element schemes 
that are stabilized by adding additional time-upwind and interface-jump terms.
These stabilizations allow him to prove a priori discretization error estimates 
in some mesh-dependent norm
as was done in the papers 
\cite{LangerNeumuellerSchafelner:2019a} or \cite{LangerSchafelner:2020a}
for linear parabolic initial-boundary value problems.

In this paper, we consider goal-oriented adaptive space-time finite element 
Galerkin discretizations of the regularized variational $p$-Laplace problem for 
{$p \in (1,\infty)$}
without any further stabilization. In goal oriented error estimations \cite{BeRa01,dolejvsi2021goal,FeiPraeZee16,MaVohYou2019,FickBrummelenZee2010,granzow2023linearization}, 
we aim 
at the estimation of
the error in 
some quantities of interest, also called goal functionals. 
For information about the treatment of multiple goal functionals at once see \cite{HaHou03,Ha08,ahuja2022multigoal,KerPruChaLaf2017,AlvParBar2013,BeEnLaWi2021a, endtmayerPAMM2022}.
In this work, we use the dual weighted residual (DWR) method \cite{BeckerRannacher1995,BeRa01},
where
the localization of the error estimator
is done by the partition-of-unity (PU) technique \cite{RiWi15_dwr}. 
For other localization techniques, we refer to \cite{BeRa01,becker2002optimal,BraackErn02}.

The DWR method additionally requires the solution of the adjoint problem 
that is linear,
	but backward in time.
So we can
basically 
use the same space-time 
finite element discretization as for the primal problem.
This is one advantage of fully unstructured space-time finite element methods
that was already used for parabolic optimal control problems 
where the adjoint problem appears in the first-order optimality system;
see, e.g., \cite{LangerSteinbachTroeltzschYang:2021b,langer2022adaptive,schafelner2022space,langer2021simultaneous}.
The goal-oriented adaptivity can now be done simultaneously in space and time 
like in the elliptic case since the time $t$ is just another variable.
The same is true for a simultaneous parallelization of the final 
numerical algorithm. 
These are two other benefits of the unstructured space-time approach 
over the more traditional time-stepping.

The remainder of the paper is organized as follows. 
In Section~\ref{Section:STVF}, we derive  the space-time variational formulation
of the regularized $p$-Laplace initial-boundary value problem \eqref{eqn:Parabolic_p-Laplace_CF}, 
and introduce some notations and preliminary results. 
This variational formulation is  the starting point for 
the space-time finite element discretization that is 
presented in Section~\ref{Section:STFEM} together with 
the Newton linearization of the discrete problem.
Section~\ref{Section:STAFEM} is devoted to the description of
the goal-oriented space-time adaptive procedure that is driven by the 
DWR method and their localization by the PU technique.
In Section~\ref{Section:NumericalResults}, we present and discuss
some numerical results for two different goal functionals.

Finally, we draw some conclusions, and give an outlook
on some further research topics.

\section{Space-Time Variational Formulations}
\label{Section:STVF}

Multiplying the parabolic $p$-Laplace equation by a test function
$v \in V := L_p((0,T); \mathring{W}_p^1(\Omega))$,
integrating over $Q$, and integrating by parts in the nonlinear 
elliptic term,
we arrive at 
the 
following variational formulation:
Find $u \in U := \{v \in V: \partial_t v \in V^*, v=0 \;\mbox{on}\; \Sigma_0\}$
such that
\begin{equation}
 \label{eqn:Parabolic_p-Laplace_VF}
\langle \partial_t u,v \rangle + ((|\nabla_{x}u|^2 + \varepsilon^2)^{\frac{p-2}{2}}\nabla_{x}u,\nabla_{x}v) = \langle f,v\rangle
\quad \forall \, v \in V,
\end{equation}
where $\langle \cdot, \cdot \rangle: V^* \times V \rightarrow \mathbb{R}$ and 
$( \cdot, \cdot): H \times H \rightarrow \mathbb{R} $ denote the duality and the $H=L_2(Q)$-inner product, respectively.
The given source function $f$ should belong to $V^* := L_q((0,T); W_q^{-1}(\Omega))$ 
that is the dual space of $V$, where $q = p/(p-1)$ denotes H\"older's conjugate of $p$.
Here and throughout the paper, we use the usual notations for Sobolev and Bochner spaces;
see, e.g., \cite{Zeidler:1990a}, 
where the norms in the trial space $U$ and the test space $V$ are defined by
\begin{equation*}
\|u\|_U = \|u\|_V + \|\partial_t u\|_{V^*}
\quad \mbox{and} \quad
\|u\|_V = \| \nabla_x u \|_{L_p(Q)}  :=  \| |\nabla_x u| \|_{L_p(Q)}  
= \Big( \int_Q \big(\sum_{i=1}^d |\partial_i v|^2\big)^{p/2}\Big)^{1/p},
\end{equation*}
respectively, with $\| f \|_{V^*} = \sup_{v\in V} \langle f, v \rangle / \|v\|_V$.
We mention that, instead of the standard $L_p$ norm for vector functions, we use the 
equivalent norm $\| |\cdot| \|_{L_p(Q)}$ that is better suited for investigating the $p$-Laplace problem.
The variational problem \eqref{eqn:Parabolic_p-Laplace_VF}  
has a unique solution $u$ that belongs to $U \cap C([0,T]; L_2(\Omega))$.
This solvability result 
for \eqref{eqn:Parabolic_p-Laplace_VF} including the case $\varepsilon=0$
follows from standard monotonicity arguments; 
see, e.g., 
\cite{Lions:1969a},
\cite{Zeidler:1990a}
and \cite{Roubicek:2013a}.
In particular,
it has recently been proved in \cite{CianchiMazya:2020a} that,
under additional regularity assumptions imposed on $\Omega$ and under the assumptions that 
\begin{equation*} 
 \label{eqn:RegularityAssumptions}
 f \in L_2(Q) \quad \mbox{and} \quad u_0 \in  \mathring{W}_p^1(\Omega),
\end{equation*}
there is a unique approximable solution to \eqref{eqn:Parabolic_p-Laplace_CF} 
for $\varepsilon = 0$
such that 
$u \in L_\infty((0,T); \mathring{W}_p^1(\Omega))$,
$\partial_t u \in L_2(Q)$ and 
$ |\nabla_x u|^{p-2} \nabla_x u \in L_2((0,T); W_2^1(\Omega))$,
i.e. the $p$-Laplace equation \eqref{eqn:Parabolic_p-Laplace_CF} 
holds in the strong sense in $L_2(Q)$. Moreover, 
the corresponding a priori estimates are proven 
in same paper; see Theorem~2.1 and Theorem~2.2 in \cite{CianchiMazya:2020a}.
Further regularity results can be found, e.g.,  in 
\cite{Zeidler:1990a}
and \cite{Roubicek:2013a}.

The nonlinear variational problem \eqref{eqn:Parabolic_p-Laplace_VF} can be rewritten 
as nonlinear operator equation: Find $u \in U$ such that
\begin{equation}
 \label{eqn:Parabolic_p-Laplace_Au=0}
    \mathcal{A}(u) := \partial_t u + A(u) - f = 0 \quad \mbox{in}\, V^*,
\end{equation}
where the nonlinear operator $\mathcal{A}: U \mapsto V^*$ is defined by the variational identity
\begin{equation*} 
\label{eqn:mathcalA}
\mathcal{A}(u)(v) = \langle \mathcal{A}(u), v \rangle
:=\langle \partial_t u,v \rangle + \langle A(u),v\rangle - \langle f,v\rangle
\quad \forall \, v \in V, \; \forall \, u \in U,
\end{equation*}
with $\langle A(u),v\rangle := ((|\nabla_{x}u|^2 + \varepsilon^2)^{\frac{p-2}{2}}\nabla_{x}u,\nabla_{x}v)$.

\section{Space-Time Finite Element Discretization and Linearization}
\label{Section:STFEM}
We are now going to construct finite element schemes as Galerkin approximations 
to \eqref{Section:STVF}
on fully unstructured simplicial decompositions of the space-time cylinder $Q$.
Let $\mathcal{T}_h = \{\Delta\}$ be a decomposition (mesh) of $Q$ into non-overlapping 
shape-regular simplicial finite elements $\Delta$ 
(triangles, tetrahedra  and pentatops for $d=1,2$, and $3$, respectively)  
such that $\overline{Q} = \bigcup_{\Delta \in \mathcal{T}_h} \overline{\Delta}$. 
{So, for simplicity, we here and in the following assume that the spacial domain $\Omega$ is polytopic.}
The discretization parameter $h$ 
should indicate that, in fact, we consider a family of finer and finer meshes, 
where the refinement can be done adaptively; see, e.g.,
\cite{BrennerScott:2008a,ErnGuermond:2004a,Steinbach:2008a}   for a precise description 
of such families of shape-regular meshes and their construction.\\
After choosing the space-time finite element mesh $\mathcal{T}_h$, we can define 
space-time finite element subspaces
\begin{equation*}
U_h = V_h := \{ v_h \in S_h^k(\overline{Q}):\, v_h=0\;\mbox{on}\;{\overline \Sigma} { \cup {\overline \Sigma}_0 } \} \subset U \subset V, 
\end{equation*}
where 
$S_h^k(\overline{Q}) = \{ v_h\in C(\overline{Q})  : v_h({\tt x}_\Delta(\cdot)) \in\mathbb{P}_k(\widehat{\Delta}),\,\forall  \Delta \in \mathcal{T}_h\}$
is nothing but the standard finite element space based on polynomials of the degree 
$k \in \mathbb{N}:= \{1,2,\dots\}$.
The regular
map ${\tt x}_\Delta(\cdot) = (x_1(\cdot),\ldots,x_{d+1}(\cdot)): \widehat{\Delta} \rightarrow \Delta$ 
maps the reference element $\widehat{\Delta}$ (unit simplex) to the finite element 
$\Delta \in \mathcal{T}_h$, and 
$\mathbb{P}_k(\widehat{\Delta})$ denotes the space of polynomials of degree $k$
on the reference element $\widehat{\Delta}$.
We look at the time $t$ as just another coordinate direction $x_{d+1}$.
Here and in our numerical experiments in Section~\ref{Section:NumericalResults},
we only consider affine-linear mappings ${\tt x}_\Delta(\cdot)$ 
{since we assumed polytopic spacial domains for simplicity.}
In this case, the finite element space $S_h^k(\overline{Q})$ 
consists of all continuous and piecewise polynomial functions.
It is clear that we can use more general curved elements realized by the mapping 
that would not only admit curved spatial domains but also changing spatial domains in time 
leading to curved, but geometrically fixed space-time ``cylinders'' $Q$.
If the Jacobians of the mapping ${\tt x}_\Delta(\cdot)$ fulfil the standard 
regularity properties, then the 
usual
approximations properties 
hold; see. e.g, \cite{BrennerScott:2008a,ErnGuermond:2004a,Steinbach:2008a}.

Once the finite element trial and test spaces $U_h$ and $V_h=U_h$ are defined,
we can look for 
a
space-time finite element solution $u_h \in U_h$ 
such that 
\begin{equation}
\label{eqn:FEScheme}
\mathcal{A}(u_h)(v_h) 
\equiv
\langle \mathcal{A}(u_h), v_h \rangle
:=\langle \partial_t u_h,v_h \rangle + \langle A(u_h),v_h\rangle - \langle f,v_h\rangle
= 0
\quad \forall \, v_h \in V_h = U_h.
\end{equation}
The space-time finite element Galerkin scheme has a unique solution. The existence follows, e.g., from
the proof of Lemma~8.95 in \cite{Roubicek:2013a} where such kind of space-time Galerkin schemes 
were used to proof the existence to weak solutions of non-linear parabolic initial-boundary 
value problems like \eqref{eqn:Parabolic_p-Laplace_CF}. The uniqueness of the finite element 
solution is a consequence of fact that the Jacobian is uniformly positive definite 
for fixed $\varepsilon$, $p$, $d$, and $h$; see 
the
estimates
given below.

Let the  finite element spaces $U_h=V_h = \mbox{span}\{\varphi_j: j=1,\ldots,N_h\}$ be spanned by 
the standard nodal finite element basis. Then the finite element solution
$u_h = \sum_{j=1}^{N_h} u_j \varphi_j$
of \eqref{eqn:FEScheme} 
can be defined by the solution $\mathbf{u}_h = (u_1,\ldots, u_{N_h})^\top{} \in \mathbb{R}^{N_h}$
of the non-linear system of finite element equations
\begin{equation}
\label{eqn:nonlinearFESystem}
\mathcal{A}_h(\mathbf{u}_h) 
= 0 \quad
\mbox{in}\; \mathbb{R}^{N_h},
\end{equation}
where the image $\mathcal{A}_h(\mathbf{u}_h)$ of the non-linear map $\mathcal{A}_h: \mathbb{R}^{N_h} \rightarrow \mathbb{R}^{N_h}$ at some given vector  $\mathbf{u}_h \in\mathbb{R}^{N_h}$ 
can be computed by the formula 
\begin{equation}
\label{eqn:A_h(u_h)}
 \mathcal{A}_h(\mathbf{u}_h) = T_h \mathbf{u}_h + \widehat{A}_h(\mathbf{u}_h) \mathbf{u}_h- \mathbf{f}_h.
\end{equation}
The $N_h \times N_h$ matrices $T_h$ and $\widehat{A}_h(\mathbf{u}_h)$ are defined by the 
variational identities
\begin{equation}
\label{eqn:matrices}
(T_h \mathbf{w}_h,\mathbf{v}_h)_{\ell_2} =  \langle \partial_t w_h,v_h \rangle \; \mbox{and} \;
(\widehat{A}_h(\mathbf{u}_h) \mathbf{w}_h,\mathbf{v}_h)_{\ell_2} = 
((|\nabla_{x}u_h|^2 + \varepsilon^2)^{\frac{p-2}{2}}\nabla_{x}w_h,\nabla_{x}v_h)\;
\end{equation}
for all $w_h, v_h \in U_h=V_h$,
respectively, 
whereas the vector
$\mathbf{f}_h = ( \langle f,\varphi_i \rangle )_{i=1,\ldots,N_h} \in \mathbb{R}^{N_h}$
can easily be computed from the given source term $f \in V^*$.
The matrix $T_h$ is non-symmetric, but non-negative, 
whereas the matrix $\widehat{A}_h(\mathbf{u}_h)$ is always symmetric and positive definite (spd)
provided that $\varepsilon^2 > 0$, 
otherwise it is always non-negative.
The unique  solvability of the nonlinear system  \eqref{eqn:A_h(u_h)}  
follows from the unique  solvability of the finite element scheme \eqref{eqn:FEScheme}.

The non-linear system \eqref{eqn:nonlinearFESystem} can be solved by means of the 
Newton method that reads as follows: For given initial guess $\mathbf{u}_h^0 \in\mathbb{R}^{N_h}$,
find $\mathbf{u}_h^{n+1} =  \mathbf{u}_h^{n} + \mathbf{w}_h^{n+1}$ 
via the Newton correction $\mathbf{w}_h^{n+1}\in\mathbb{R}^{N_h}$ that is nothing but the unique 
solution of the linear system
\begin{equation*}
\label{eqn:Newton_h}
\mathcal{A}'_h(\mathbf{u}_h^{n}) \mathbf{w}_h^{n+1} = - \mathcal{A}_h(\mathbf{u}_h^{n}), \quad n=0,1,\ldots,
\end{equation*}
with the Jacobian 
$\mathcal{A}'_h(\mathbf{u}_h^{n}) = T_h + \widehat{A}_h(\mathbf{u}_h^{n}) + \widehat{A}'_h(\mathbf{u}_h^{n})$,
where the matrices $T_h$ and $\widehat{A}_h(\mathbf{u}_h^{n})$ are given by \eqref{eqn:matrices},
and the matrix $\widehat{A}'_h(\mathbf{u}_h^{n})$ is defined by the identity

\begin{equation*}
\label{eqn:matrixA'}
(\widehat{A}'_h(\mathbf{u}_h^{n}) \mathbf{w}_h,\mathbf{v}_h)_{\ell_2} =  
(p-2) (((|\nabla_{x}{u}_h^{n}|^2 + \varepsilon^2)^{\frac{p-4}{2}} \nabla_x {u}_h^{n} (\nabla_x {u}_h^{n})^\top{}) \nabla_{x} w_h,\nabla_{x}v_h)
\end{equation*}
for all $w_h, v_h \in U_h$ associated to the vectors $\mathbf{w}_h,\mathbf{v}_h \in \mathbb{R}^{N_h}$.
The Jacobian is positive definite 
for $\varepsilon>0$.
Indeed, 
for $p > 1$, we have
\begin{align*}
(\mathcal{A}'_h(\mathbf{u}_h^{n}) \mathbf{w}_h,\mathbf{w}_h)_{\ell_2} 
    =&\;(\partial_t w_h,w_h) +((|\nabla_{x}u_h^{n}|^2 + \varepsilon^2)^{\frac{p-2}{2}}\nabla_{x}w_h,\nabla_{x}w_h)\\
     &+ (p-2) (((|\nabla_{x}{u}_h^{n}|^2 + \varepsilon^2)^{\frac{p-4}{2}} \nabla_x {u}_h^{n} (\nabla_x {u}_h^{n})^\top{}) \nabla_{x} w_h,\nabla_{x}w_h)\\
    & \ge \frac{1}{2} \|w_h(\cdot,T)\|_{L_2(\Omega)}^2 + c(p,\varepsilon) \|\nabla_{x} w_h\|_{L_2(Q)}^2
    \ge c(p,\varepsilon) c h^{d+1} (\mathbf{w}_h,\mathbf{w}_h)_{\ell_2}
\end{align*}    
for all $w_h \in U_h$ corresponding to the coefficient vector $\mathbf{w}_h \in \mathbb{R}^{N_h}$
via the finite element isomorphism, 
where $c(p,\varepsilon) = (p-1) \varepsilon^{p-2}$  and $c(p,\varepsilon) = \varepsilon^{p-2}$ 
for $p \in (1,2)$ and $p \in [2,\infty)$, respectively.
Thus, the Newton method is always well-defined.

\section{Goal-Oriented Adaptive Space-Time  Finite Element Methods}
\label{Section:STAFEM}

In many applications, the solution $u$ itself is often not of primary interest,
but some quantity of interest represented by some (possibly nonlinear) functional
$J:U\mapsto \mathbb{R}$ evaluated at the solution $u$.
Such quantities could be 
the average of the solution or the gradient of the solution in some 
subregion of the space-time cylinder, or a regularized point evaluation; 
see also the quantity of interest
used in our numerical experiments in Section~\ref{Section:NumericalResults}.
The aim is to approximate $J(u)$ 
as well as  
possible using 
as few as possible degrees of freedom.
Of course, instead of $J(u)$, we can only compute the functional at the finite 
element solution $u_h$, i.e. $J(u_h)$.
Therefore, we are interested in evaluation of the error $J(u)-J(u_h)$.
To estimate the error, we use the dual weighted residual method as described in \cite{BeckerRannacher1995,BeRa01,RanVi2013}.
Other GOAST methods using time-slabs or a recuded order basis can be found in \cite{KoeBruBau2019a,roth2022tensor,fischer4420888more}.

\subsection{The primal and the adjoint problems}
\label{Subsection:PrimalAndAjointProblems}
The primal problem 
is nothing but
the variational problem \eqref{eqn:Parabolic_p-Laplace_Au=0},
and it defines the solution $u \in U$.
Its   discretization is given by the finite element scheme \eqref{eqn:FEScheme}
defining the finite element solution $u_h \in U_h$ approximating $u$.

To connect the functional $J$ with the initial-boundary value problem \eqref{eqn:Parabolic_p-Laplace_Au=0}, 
we consider the adjoint problem. 
The {formal} adjoint problem reads as follows:
Find $z \in V$ such that
\begin{equation}
\label{eqn:AdjointProblem_VF} 
\mathcal{A}'(u)(v,z) := \langle\partial_t v, z\rangle + A'(u)(v,z)= J'(u)(v) \quad \forall v \in U,
\end{equation}
where $u$ solves the primal problem \eqref{eqn:Parabolic_p-Laplace_Au=0},
 and 
\begin{equation*}
\label{eqn:AdjointOperator} 
A'(u)(v,z) =  (((|\nabla_{x}u|^2 + \varepsilon^2)^{\frac{p-2}{2}} I 
   + (p-2) (|\nabla_{x}u|^2 + \varepsilon^2)^{\frac{p-4}{2}} \nabla_x u (\nabla_x u)^\top{})
   \nabla_{x} v,\nabla_{x}z).
\end{equation*}

To obtain $z$ we still have to solve a linear 
PDE problem
	that is nothing but a well-posed backward parabolic problem.
Therefore, the adjoint problem \eqref{eqn:AdjointProblem_VF} must also be discretized 
in order to obtain a finite element approximation $z_h$ to the mesh sensitivities $z \in V$:
Find $z_h \in V_h$ such that 
\begin{equation}
\label{eqn:discrete adjoint problem}
\mathcal{A}'(u_h)(v_h,z_h) := \langle\partial_t v_h, z_h\rangle + A'(u_h)(v_h,z_h)= J'(u_h)(v_h) \quad \forall v_h \in U_h,
\end{equation}
where $u_h$ solves the discrete primal problem \eqref{eqn:FEScheme}, and $U_h = V_h$.
We note that the discrete adjoint problem \eqref{eqn:discrete adjoint problem} 
has a unique solution $z_h \in V_h$ since $\mathcal{A}'(u_h)$ is positive definite.


\subsection{An error identity}
\label{Subsection:AnErrorIdentity}

With the solutions to the primal and adjoint problems, we can show the following error identity 
for continuously differentiable operators in general:

\begin{theorem}\label{Theorem: Error Representation}
	Let us assume that $\mathcal{A} \in \mathcal{C}^3(U,V^*)$ and $J \in \mathcal{C}^3(U,\mathbb{R})$, where $V^*$ is the dual space of $V$. 
	Let $u$ be the solution of the primal problem \eqref{eqn:Parabolic_p-Laplace_Au=0}, and $z$ be the solution of the adjoint problem \eqref{eqn:AdjointProblem_VF}. 
	Then
	the error representation formula
	\begin{align} \label{Eq:Error Representation}
	\begin{split}
	J(u)-J(\tilde{u})&= \frac{1}{2}\rho(\tilde{u})(z-\tilde{z})+\frac{1}{2}\rho^*(\tilde{u},\tilde{z})(u-\tilde{u}) 
	-\rho (\tilde{u})(\tilde{z}) + \mathcal{R}^{(3)},
	\end{split}
	\end{align}
	holds for arbitrary fixed  $\tilde{u} \in U$ and $ \tilde{z} \in V$, where
	$\rho(\tilde{u})(\cdot) := -\mathcal{A}(\tilde{u})(\cdot)$,
	$\rho^*(\tilde{u},\tilde{z})(\cdot) := J'(u)-\mathcal{A}'(\tilde{u})(\cdot,\tilde{z})$,
	and the remainder term
	\begin{equation*}
	\begin{split}	\label{Eq:Error Estimator: Remainderterm}
	\mathcal{R}^{(3)}:=\frac{1}{2}\int_{0}^{1}\big[J'''(\tilde{u}+se)(e,e,e)
	-\mathcal{A}'''(\tilde{u}+se)(e,e,e,\tilde{z}+se^*)
	-3\mathcal{A}''(\tilde{u}+se)(e,e,e)\big]s(s-1)\,ds,
	\end{split} 
	\end{equation*}
	with $e=u-\tilde{u}$ and $e^* =z-\tilde{z}$.
\end{theorem}
\begin{proof}
	The detailed proof can be found in
	\cite{BeRa01,RanVi2013,EnLaWi20,EnLaWi18}.
\end{proof}
We note that we can weaken the differentiability assumptions to the corresponding 
differentiability on the line between $\tilde{u}$ and $u$.
\begin{remark}
	While $\tilde{u}$ and $\tilde{z}$ can be arbitrary functions from $U$ and $V$, we will later fix them as the finite element solutions of the primal and adjoint problem, respectively.
	Moreover, we note that we will later need the  continuous differentiability of $\mathcal{A}$
	only in the finite element subspaces; see Subsection~\ref{Subsection:ErrorEstimatorLocalization}.
\end{remark}

Of course, \eqref{Eq:Error Representation} can be used as an error estimation. 
However, we do not know the exact solutions $u$ and $z$ 
of the primal and the adjoint problems, respectively. 
Therefore, \eqref{Eq:Error Representation} is not computable.
	

We now replace
$u$ and $z$ in \eqref{Eq:Error Representation} 
by computable quantities. 
In order to accomplish this, we introduce enriched 
finite element
spaces $U_h^{(2)}$ and $V_h^{(2)}$ with the properties $U_h \subset U_h^{(2)} \subset U$ and $V_h \subset V_h^{(2)} \subset V$.
On the basis of these new spaces,
we form the enriched primal and adjoint finite element problems.
The enriched primal finite element problem 
can be formulated as follows:
Find $u_h^{(2)} \in U_h^{(2)}$ such that 
\begin{equation}
\label{eqn: enriched discrete primal problem.}
\mathcal{A}(u_h^{(2)})(v_h^{(2)}) = \langle \mathcal{A}(u_h^{(2)}), v_h^{(2)} \rangle
:=\langle \partial_t u_h^{(2)},v_h^{(2)} \rangle + \langle A(u_h^{(2)}),v_h^{(2)}\rangle - \langle f,v_h^{(2)}\rangle
\quad \forall \, v_h^{(2)} \in V_h^{(2)}.
\end{equation}
Let $u_h^{(2)}$ be the solution of \eqref{eqn: enriched discrete primal problem.}. 
Then the enriched adjoint finite element problem reads as follows: Find $z_h^{(2)} \in V_h^{(2)}$ such that 
\begin{equation}
\label{eqn: enriched discrete adjoint problem}
\mathcal{A}'(u_h^{(2)})(v_h^{(2)},z_h^{(2)}) := \langle\partial_t v_h^{(2)}, z_h^{(2)}\rangle + A'(u_h^{(2)})(v_h^{(2)},z_h^{(2)})= J'(u_h^{(2)})(v_h^{(2)}) \quad \forall v_h^{(2)} \in U_h^{(2)}.
\end{equation}
If we now replace $u$ and $z$ by our enriched finite element solutions $u_h^{(2)}$ and $z_h^{(2)}$ 
in \eqref{Eq:Error Representation}, 
we arrive at the approximate error representation
\begin{align} \label{eqn: etafull}
J(u)-J(\tilde{u})\approx \frac{1}{2}\underbrace{\rho(\tilde{u})\left(z_h^{(2)}-\tilde{z}\right)}_{=:\eta_{h,p}}+\frac{1}{2}\underbrace{\rho^*(\tilde{u},\tilde{z})\left(u_h^{(2)}-\tilde{u}\right)}_{=:\eta_{h,a}} 
-\underbrace{\rho (\tilde{u})(\tilde{z})}_{=:\eta_{k}} + \underbrace{\mathcal{R}^{(3),(2)}}_{=:\eta_{\mathcal{R},2}},
\end{align}
where 
	\begin{equation*}
\begin{split}	\label{Eq:Error Estimator: Remainderterm Enriched}
\mathcal{R}^{(3),(2)}:=\frac{1}{2}\int_{0}^{1}\big[J'''(\tilde{u}+s\hat{e})(\hat{e},\hat{e},\hat{e})
-\mathcal{A}'''(\tilde{u}+s\hat{e})(\hat{e},\hat{e},\hat{e},\tilde{z}+s\hat{e}^*)
-3\mathcal{A}''(\tilde{u}+s\hat{e})(\hat{e},\hat{e},\hat{e})\big]s(s-1)\,ds,
\end{split} 
\end{equation*}
with $\hat{e}=u_h^{(2)}-\tilde{u}$ and $\hat{e}^* =z_h^{(2)}-\tilde{z}$.
\begin{remark}
	As above $\tilde{u}$ and $\tilde{z}$ are meant to be the solutions of \eqref{eqn:FEScheme} and  \eqref{eqn:discrete adjoint problem} or close approximations to them.
\end{remark}%
This enriched approximation is also used in \cite{BeRa01,BaRa03,blum2003posteriori,HaHou03,Ha08,bruchhauser2017numerical,KoeBruBau2019a,EnWi17,EnLaNeiWoWi2020,BeEnWi2021P}. A comparison between $h$-enrichment and $p$-enrichment can be found in \cite{endtmayer2021hierarchical}.

\subsection{A different representation}
In this section, we refine the analysis of the previous section such that the Fr\'{e}chet differentiability of the operator and goal functional are just required on the enriched spaces.

\begin{assumption}[Saturation assumption for the goal functional; see \cite{EnLaWi20}]
	\label{Assumption: Saturation}
	Let $u_h^{(2)}$ be the solution of the enriched  primal problem \eqref{eqn: enriched discrete primal problem.}, and $z_h^{(2)}$ be the solution of the enriched adjoint  problem \eqref{eqn: enriched discrete adjoint problem}. We assume that $\tilde{u}$ and $\tilde{z}$ are the solutions of \eqref{eqn:FEScheme} and  \eqref{eqn:discrete adjoint problem} 
	or close approximations to them. 
	Then there is a constant $b<1$ such that 
	$$|J(u)-J(u_h^{(2)})| < b \, |J(u)-J(\tilde{u})|.$$ 
\end{assumption}

\begin{corollary}\label{Corollary: Error Representation}
Let $\mathcal{A}: U \mapsto V^*$ and $J:U\mapsto \mathbb{R}$.
Furthermore, let us assume that $J(u) \in \mathbb{R}$, where $u\in U$ solves the primal problem \eqref{eqn:Parabolic_p-Laplace_Au=0}.
Additionally, let $u_h^{(2)}$ be the solution of the enriched  primal problem \eqref{eqn: enriched discrete primal problem.}, and $z_h^{(2)}$ be the solution of the enriched adjoint  problem \eqref{eqn: enriched discrete adjoint problem}. 
Moreover, let  $\mathcal{A}_h \in \mathcal{C}^3(U_h^{(2)},V_h^{(2)*})$ and $J_h \in \mathcal{C}^3(U_h^{(2)},\mathbb{R})$ such that 
for all $u_h^{(2)},\psi_h^{(2)} \in U_h^{(2)}$ and $\phi_h^{(2)} \in V_h^{(2)}$, 
the equalities 
\begin{align}
	\label{eqn: A=Ah}
	\mathcal{A}(u_h^{(2)})(\phi_h^{(2)})&=\mathcal{A}_h(u_h^{(2)})(\phi_h^{(2)}), \\
	\label{eqn: A'=Ah'}
	\mathcal{A}'(u_h^{(2)})(\psi_h^{(2)},\phi_h^{(2)})&=\mathcal{A}'_h(u_h^{(2)})(\psi_h^{(2)},\phi_h^{(2)}),\\
	\label{eqn: J=Jh}
	J(\psi_h^{(2)})&=J_h(\psi_h^{(2)}),\\
	\label{eqn: J'=Jh'}
	J'(\psi_h^{(2)})&=J'_h(\psi_h^{(2)}),
\end{align}
 are fulfilled. 
Here, $V_h^{(2)*}$ denotes the dual space of $V_h^{(2)}$. 
Then
the error representation formula
\begin{align*} \label{Eq:Error RepresentationEnriched}
\begin{split}
J(u)-J(\tilde{u})&= J(u)-J(u_h^{(2)})+ \frac{1}{2}\rho(\tilde{u})(z_h^{(2)}-\tilde{z})+\frac{1}{2}\rho^*(\tilde{u},\tilde{z})(u_h^{(2)}-\tilde{u}) 
-\rho (\tilde{u})(\tilde{z}) + \mathcal{R}_h^{(3)},
\end{split}
\end{align*}
holds for arbitrary fixed  $\tilde{u} \in U_h^{(2)}$ and $ \tilde{z} \in V_h^{(2)}$, where
$\rho(\tilde{u})(\cdot) := -\mathcal{A}(\tilde{u})(\cdot)$ and
$\rho^*(\tilde{u},\tilde{z})(\cdot) := J'(u)-\mathcal{A}'(\tilde{u})(\cdot,\tilde{z})$.
\end{corollary}
\begin{proof}
Since $u_h^{(2)}$ solves the enriched problem \eqref{eqn: enriched discrete primal problem.}, we know 
\begin{equation*}
	\mathcal{A}(u_h^{(2)})(v_h^{(2)}) =0  \quad \forall v_h^{(2)} \in V_h^{(2)},
\end{equation*}
and since $z_h^{(2)}$ solves the enriched adjoint  problem \eqref{eqn: enriched discrete adjoint problem} we have 
\begin{equation*}
	\mathcal{A}'(u_h^{(2)})(v_h^{(2)},z_h^{(2)})= J'(u_h^{(2)})(v_h^{(2)}) \quad \forall v_h^{(2)} \in U_h^{(2)}.
\end{equation*}
If we exploit the equalities \eqref{eqn: A=Ah},\eqref{eqn: A'=Ah'} and \eqref{eqn: J'=Jh'}, we observe that $u_h^{(2)}$ and $z_h^{(2)}$ fulfill
\begin{equation*}
	\mathcal{A}_h(u_h^{(2)})(v_h^{(2)}) =0  \quad \forall v_h^{(2)} \in V_h^{(2)},
\end{equation*}
and
\begin{equation*}
	\mathcal{A}'_h(u_h^{(2)})(v_h^{(2)},z_h^{(2)})= J'_h(u_h^{(2)})(v_h^{(2)}) \quad \forall v_h^{(2)} \in U_h^{(2)}.
\end{equation*}
This allows us to apply Theorem \ref{Theorem: Error Representation} with $\mathcal{A}=\mathcal{A}_h$ and $J=J_h$ where we obtain 
\begin{equation}
	\label{eqn: pure discrete result}
		J_h(u_h^{(2)})-J_h(\tilde{u})= \frac{1}{2}\rho_h(\tilde{u})(z_h^{(2)}-\tilde{z})+\frac{1}{2}\rho^*_h(\tilde{u},\tilde{z})(u_h^{(2)}-\tilde{u}) 
	-\rho_h (\tilde{u})(\tilde{z}) + \mathcal{R}_h^{(3)},
\end{equation}
with $\rho_h(\tilde{u})(\cdot) := -\mathcal{A}_h(\tilde{u})(\cdot)$, 
$\rho^*_h(\tilde{u},\tilde{z})(\cdot) := J'_h(u)-\mathcal{A}'_h(\tilde{u})(\cdot,\tilde{z})$, and 
\begin{equation*}
	\mathcal{R}^{(3)}_h:=\frac{1}{2}\int_{0}^{1}\big[J_h'''(\tilde{u}+s\hat{e})(\hat{e},\hat{e},\hat{e})
	-\mathcal{A}_h'''(\tilde{u}+s\hat{e})(\hat{e},\hat{e},\hat{e},\tilde{z}+s\hat{e}^*)
	-3\mathcal{A}_h''(\tilde{u}+s\hat{e})(\hat{e},\hat{e},\hat{e})\big]s(s-1)\,ds.
\end{equation*}
Here, $\hat{e}=u_h^{(2)}-\tilde{u}$ and $\hat{e}^* =z_h^{(2)}-\tilde{z}$.
Using \eqref{eqn: A=Ah},\eqref{eqn: A'=Ah'}and \eqref{eqn: J'=Jh'}, one can observe that $\rho_h=\rho$ and $\rho_h^*=\rho^*$ on the enriched spaces. 
In combination with \eqref{eqn: pure discrete result} and \eqref{eqn: J=Jh},
this leads to that
\begin{align*}
	J_h(u_h^{(2)})-J_h(\tilde{u})= \frac{1}{2}\rho_h(\tilde{u})(z_h^{(2)}-\tilde{z})+\frac{1}{2}\rho^*_h(\tilde{u},\tilde{z})(u_h^{(2)}-\tilde{u}) 
	-\rho_h (\tilde{u})(\tilde{z}) + \mathcal{R}_h^{(3)},\\
\end{align*}
is equivalent to
\begin{align*}
	J(u)-	J(u_h^{(2)})+\frac{1}{2}\rho_h(\tilde{u})(z_h^{(2)}-\tilde{z})+\frac{1}{2}\rho^*_h(\tilde{u},\tilde{z})(u_h^{(2)}-\tilde{u} 
	-\rho_h (\tilde{u})(\tilde{z}) + \mathcal{R}_h^{(3)}&=\\
	J(u)-	J(u_h^{(2)})+	J_h(u_h^{(2)})-J_h(\tilde{u})&=J(u)-J(\tilde{u}).
\end{align*}
This concludes the proof.
\end{proof}
\begin{remark}
	In contrast to Theorem \ref{Theorem: Error Representation}, 
	Corollary \ref{Corollary: Error Representation} does not require differentiability of $\mathcal{A}$ and $J$ in the Sobolev spaces $U,V$, but in the discrete spaces $U_h^{(2)}, V_h^{(2)}.$ Furthermore, we observe that 
	$$J(u)-J(\tilde{u})= J(u)-J(u_h^{(2)})+\frac{1}{2}(\eta_{h,p}+\eta_{h,a})+\eta_{k}+\eta_{\mathcal{R}},$$
	with $\eta_{\mathcal{R}}:=\mathcal{R}_h^{(3)}$.
If the saturation assumption is fulfilled, we have 
$$(1-b)\,|J(u)-J(\tilde{u})| \leq |\frac{1}{2}(\eta_{h,p}+\eta_{h,a})+\eta_{k}+\eta_{\mathcal{R}}|.$$
\end{remark}
\subsection{The discretization error estimator and its localization}
\label{Subsection:ErrorEstimatorLocalization}
In this section, we discuss the different parts of the error estimator.

The discretization error estimator $\eta_{h}$ is given by 
\begin{equation*} \label{eqn: etah}
\eta_{h}:=\frac{1}{2}\left(\eta_{h,p}+\eta_{h,a}\right),
\end{equation*}
where $\eta_{h,p}$ defined in \eqref{eqn: etafull} is the primal part of the error estimator, 
and $\eta_{h,a}$ given in \eqref{eqn: etafull} is the adjoint part of the error estimator. 
In the literature, 
often only the primal part is used.
Here, we use both parts, since it  was shown in \cite{EnLaWi18} 
that, especially for stationary $p$-Laplace problems, both terms are beneficial. 
As suggested in \cite{RanVi2013,endtmayer2020multi}, $\eta_{h}$ represents the discretization error.
Therefore, 
we will use it
to drive the mesh adaptation process. 
For this, $\eta_{h}$ must be localized. 
Here, we use the partition of unity technique 
proposed
in \cite{RiWi15_dwr}. Let $\{ \Psi_i \}_{i=1}^{N}$ be a set of functions with the property $\sum_{i=1}^{N} \Psi_{i}\equiv 1$. Then we have $\eta_{h}=\sum_{i=1}^{N} \eta_{i}$ where 
\begin{equation} \label{eta_i}
\eta_{i}:=\frac{1}{2}{\rho(\tilde{u})\left((z_h^{(2)}-\tilde{z})\Psi_i\right)}+\frac{1}{2}{\rho^*(\tilde{u},\tilde{z})\left((u_h^{(2)}-\tilde{u})\Psi_i\right)}.
\end{equation} 
In our numerical experiments, we choose $\{ \Psi_i \}_{i=1}^{N}$ to be the basis functions of the lowest-order discontinuous finite space, i.e.\ piecewise constant functions.
The local error estimators $\eta_{i}$ are then used for mesh adaption.

The iteration error estimator $\eta_k$ represents the iteration error, as suggested in \cite{RanVi2013,EnLaWi20}. If $\tilde{u}=u_h$ is the exact solution of the discretized primal problem, then $\eta_{k}=0$. It can be used to stop the nonlinear solver as done in \cite{RanVi2013,EnLaWi20}.

The part $\eta_{\mathcal{R}}$ is of higher order. 
Therefore, it is usually neglected in the literature. 
For the regularized stationary $p$-Laplace problem, this part was numerically analyzed 
in \cite{EnLaWi20}. Indeed, the results showed that $\eta_{\mathcal{R}}$ 
can be neglected.


\subsection{The final adaptive algorithm}
\label{Subsection:AdaptiveAlgorithm}

The DWR driven, goal-oriented space-time adaptive finite element approach described above
can be summarised in form of 
Algorithm~\ref{Full Adaptive Algorithm}.

\begin{algorithm}
\caption{The adaptive space-time algorithm \label{Full Adaptive Algorithm}}
\begin{algorithmic}[1] 
\Repeat
\State solve the nonlinear primal problem as in \eqref{eqn:FEScheme} using some nonlinear solver, \label{Alg:primal}
\State solve the linear adjoint problem as in \eqref{eqn:discrete adjoint problem}  using some linear solver, \label{Alg:adjoint}
\State solve the nonlinear  enriched primal problem as in \eqref{eqn: enriched discrete primal problem.} using some nonlinear solver, \label{Alg:primal+}
\State solve the linear enriched adjoint problem as in \eqref{eqn: enriched discrete adjoint problem} using some linear solver, \label{Alg:adjoint+}
\State compute the elementwise contributions  via PU-technique as in \eqref{eta_i}; see \cite{RiWi15_dwr},
\State select a set of marked elements $M$ using some marking strategy; e.g., D\"orfler marking \cite{Doerfler:1996a}, 
\State $\mathcal{T}_{k+1} \gets \Call{Refine}{\mathcal{T}_k, M} $,
\State $ k \gets k + 1$,
\Until some stopping criterion is fulfilled.
\end{algorithmic}
\end{algorithm}

The numerical results presented in Section~\ref{Section:NumericalResults} are based 
on an implementation of this algorithm. There we give some more specific information 
on the practical realization of the algorithm, in particular, on the nonlinear and 
linear solvers used.

The Algorithm~\ref{Full Adaptive Algorithm} requires the solution of two non-linear 
and two linear systems of finite element equations. This seems to be quite expensive.
However, the goal-oriented adaptive approach in connection with a nested iteration 
setting can considerable reduce the cost in comparison with a naive approach. 
Moreover, the solution of the enriched systems can be avoided as discussed in 
Subsection~\ref{Subsection:AnErrorIdentity}.
In our numerical experiments presented in the next section, we stop the adaptive
process as soon as we reach a total of \num{e6} dofs.

%
\section{Numerical Results}
\label{Section:NumericalResults}

We implemented the adaptive finite element method described in the previous section 
in our space-time finite element code 
which is based on the finite element library MFEM \cite{mfem:2021,mfem-web}. 
The linear solvers are either provided by the solver library \emph{hypre}\footnote{\url{https://computing.llnl.gov/projects/hypre-scalable-linear-solvers-multigrid-methods}} or by the sparse direct solver MUMPS\footnote{\url{https://mumps-solver.org/index.php}} 
\cite{AmestoyDuffLExcellentKoster:2001a}, 
which is used via the library PETSc\footnote{\url{https://petsc.org/release/}} 
\cite{petsc-user-ref}.
To measure the quality of our error estimates, 
we use so-called 
efficiency
indices $\mathrm{I_{eff}}$ for 
the full discretization error estimator $\eta_{h}$, 
the adjoint part $\eta_{h,a}$, and 
the primal part $\eta_{h,p}$, which are defined by the formulas
	\begin{equation*}
	\mathrm{I}_{\mathrm{eff},h}:=\frac{\eta_{h}}{J(u)-J(u_h)}, \quad \mathrm{I}_{\mathrm{eff},a}:=\frac{\eta_{h,a}}{J(u)-J(u_h)}, 
	\quad \mbox{and} \quad 
	\mathrm{I}_{\mathrm{eff},p}:=\frac{\eta_{h,p}}{J(u)-J(u_h)}. 
	\end{equation*}

We mention that we use different mesh refinement techniques for $d=2$ and $d=3$. 
{More precisely, we apply octasection for the refinement of tetrahedra ($d+1=3$) \cite{mfem:2021},
and the newest vertex bisection for the refinement of pentatopes ($d+1=4)$ \cite{Stevenson:2008a}.
}

As already described above,  we solve the nonlinear finite element equations by means of Newton's method.
We start with a damped version where the damping parameter is chosen by a simple, residual based line-search.
Close to the solution, the line search will always accept a damping parameter equal to $1$, 
and the damped version turns into the pure Newton method. 
When comparing Newton's method with different interior solvers for the Jacobian systems arising at each Newton step,
we always use the same pseudo-random initial guess for the Newton iteration.
However, in practice when one wants to solve the nonlinear systems efficiently, 
we will interpolate the computed solution from the current mesh to its adaptively refined mesh.
This interpolated solution will then serve as initial guess for the Newton solver on the next mesh level.
This nested iteration setting considerably speeds up the overall solution process.

The linear systems arising in Newton's method are  solved either by means of the Generalized Minimal Residual (GMRES) method \cite{Saad:2003a}, or by a sparse direct solver \cite{Davis:2006a}. 
We do not use any restarts for GMRES, and stop either when the initial residual has been reduced by a factor of $ 10^{-8}$, or after 100 iterations.
This procedure can certainly be improved by adapting the accuracy of the inner GMRES iteration 
to the error reduction of the Newton iteration 
\cite{Deuflhard2011}
or the error reduction in $\eta_{k}$ \cite{RanVi2013}.

\subsection{Convergence Studies for a Smooth Solution}
\label{ex:regular, smooth}

We now consider the regularized parabolic $p$-Laplace equation 
\eqref{eqn:Parabolic_p-Laplace_CF}
with the particular choices $ p = 4 $, 
$ \varepsilon \in \{1, 10^{-5}, 10^{-10}\} $, 
$Q = (0,1)^{d+1}$, and the manufactured smooth solution
\begin{equation}
\label{eqn:Example5.1:ExactSolution}
u(x,t) = t^2\,e^{t}\,\prod_{i=1}^d \sin(x_i\,\pi). 
\end{equation}
The source term $f$ is computed accordingly. 
Before we numerically study the goal-oriented space-time adaptivity proposed 
in Section~\ref{Section:STAFEM}, we examine the performance of our space-time finite element method 
with respect to (wrt)
uniform $h$-refinements as well as 
wrt 
different polynomial degrees $k$ of 
the finite element shape functions used.
In Figure~\ref{fig:regular, smooth:convergence}, we present the overall convergence history of the discretization error measured in the  $L_2(Q)$-norm as well as in the
$L_2(0,T;W^{1,2}(\Omega)) $-norm
for linear ($k=1$) and quadratic ($k=2$) shape functions in the case $d=2$.
We observe that the optimal rate provided by the approximation power of the finite element spaces
is always achieved 
for the $L_2(0,T;W^{1,2}(\Omega)) $-norm whereas the convergence order of the $L_2(Q)$-norm is reduced by one for quadratic shape functions.
For $d=3$, 
the error measured in the $L_2(Q)$- and $L_2(0,T;W^{1,2}(\Omega))$-norm decays with the optimal rate wrt the corresponding finite element spaces; see Figure~\ref{fig:regular, smooth:convergence:4d}.
\begin{figure}[!htb]
	\centering
	\includegraphics[width=.95\linewidth]{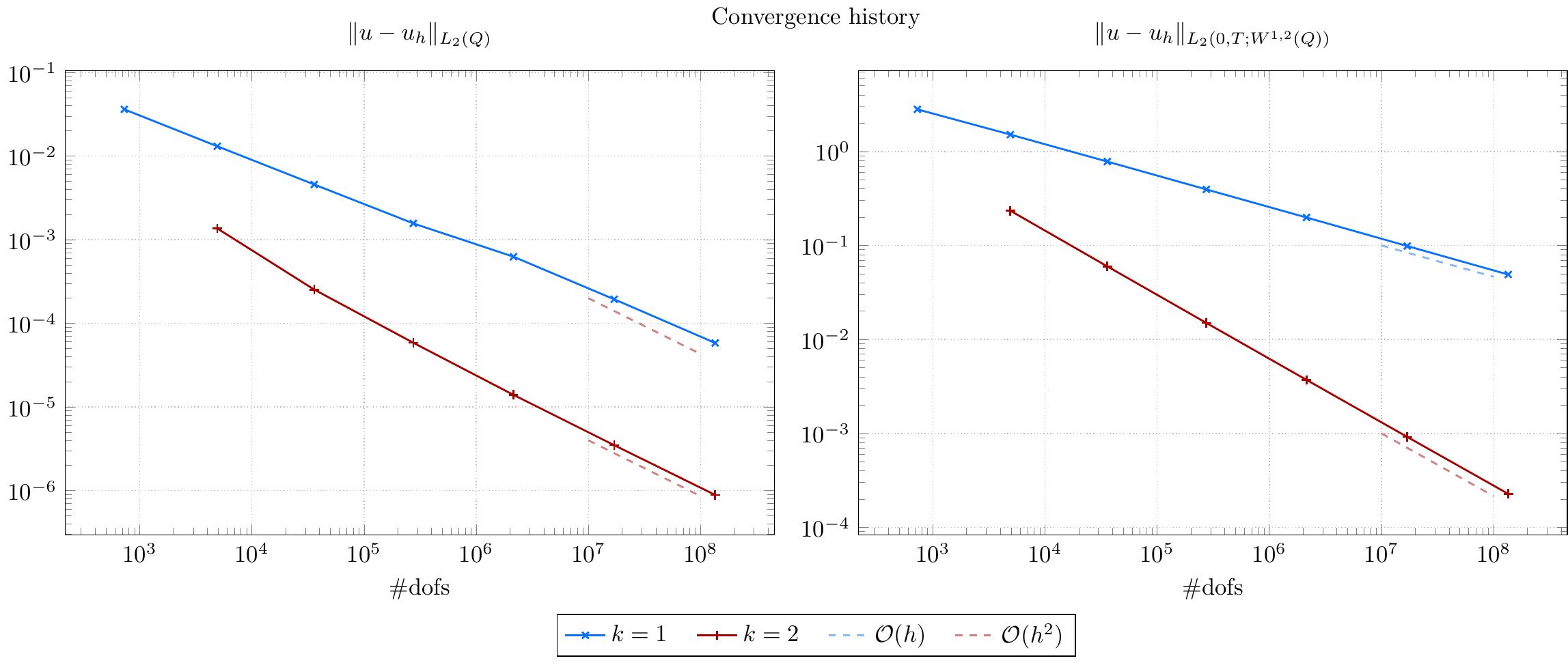}
	\caption{Example~\ref{ex:regular, smooth} ($d=2$): 
	Convergence history of the discretization errors in different norms for linear ($k=1$) and quadratic ($k=2$) shape functions.
	}
	\label{fig:regular, smooth:convergence}
\end{figure}

\begin{figure}[!htb]
	\centering
	\includegraphics[width=.95\linewidth]{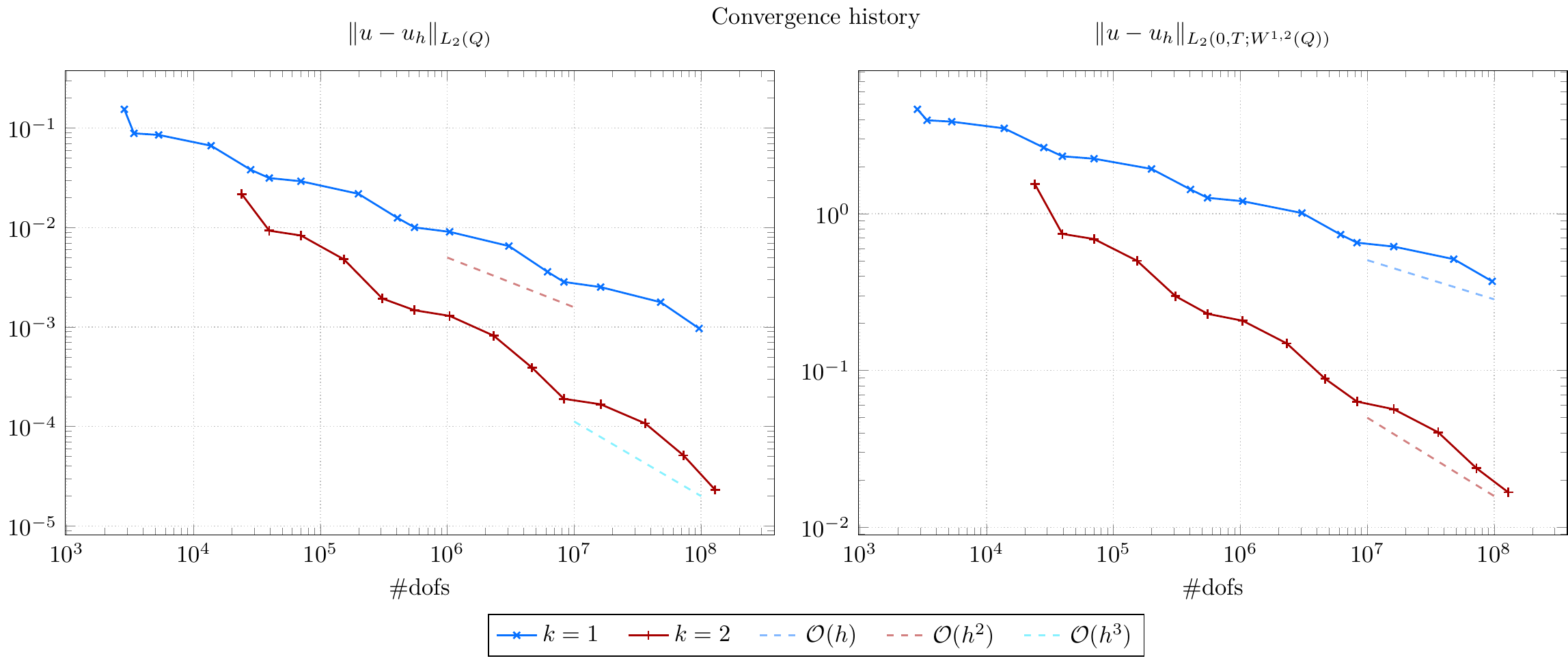}
	\caption{Example~\ref{ex:regular, smooth} ($d=3$): 
	Convergence history of the discretization errors in different norms for linear ($k=1$) and quadratic ($k=2$) 
	shape functions.
	}
	\label{fig:regular, smooth:convergence:4d}
\end{figure}
Finally, 
we study
the 
performance of the nested Newton solver. As already mentioned earlier, 
we can solve the linear systems arising in Newton's method by means of a sparse direct solver, 
or by means of an iterative method, e.g., 
the preconditioned GMRES method. 
The preconditioner is constructed by algebraic multigrid (AMG),
where we use a standard V-cycle with one pre- and one post-smoothing step;  
cf. \cite{Hackbusch:1985a,TrottenbergOosterleeSchueller:2001a}. 
Moreover, we will 
numerically investigate
the influence of the regularization parameter $ \varepsilon $ on the 
solution process.
In Table~\ref{tab:regular, smooth:h-scaling}, we present scaling results wrt uniform mesh refinement. 
Here, we can observe that the regularization parameter has only a mild influence on Newton's method, but does affect the convergence behaviour of the inner GMRES solver;  
compare, e.g., the total number of GMRES iterations in brackets 
for refinement level $ \ell = 3 $ in Table~\ref{tab:regular, smooth:h-scaling}.
\begin{table}[!ht]
	\centering%
	\caption{Example~\ref{ex:regular, smooth}: Scaling of the (damped) Newton solver with total number of inner solves in brackets, for $d+1=3$ and $k=1$;
	 using 256 cores of HPC Cluster \texttt{RADON1}\protect\footnotemark.}
	\label{tab:regular, smooth:h-scaling}
	\begin{tabular}{clrlcrlcrlc}
		\toprule \multirow{2}{*}{$ \ell $} & \multirow{2}{*}{$ \#\mathrm{dofs} $} & \multicolumn{3}{c}{$ \varepsilon = 1$} & \multicolumn{3}{c}{$ \varepsilon = 10^{-5}$} & \multicolumn{3}{c}{$ \varepsilon = 10^{-10}$}                                                                            \\
		                                   &                                      & \multicolumn{2}{c}{\textbf{GMRES}}     & {\textbf{MUMPS}}                             & \multicolumn{2}{c}{\textbf{GMRES}}            & {\textbf{MUMPS}} & \multicolumn{2}{c}{\textbf{GMRES}} & {\textbf{MUMPS}} \\ \midrule %
		0                                  & \num{4913}                           & \num{11} & \hspace*{-.85em}(\num{322}) & \num{11}                                     & \num{12} & \hspace*{-.85em}(\num{370})        & \num{12}         & \num{12} & \hspace*{-.85em}(\num{370}) & \num{12} \\
		1                                  & \num{35937}                          & \num{3}  & \hspace*{-.85em}(\num{70})  & \num{3}                                      & \num{3}  & \hspace*{-.85em}(\num{185})        & \num{3}          & \num{3}  & \hspace*{-.85em}(\num{185}) & \num{3}  \\
		2                                  & \num{274625}                         & \num{2}  & \hspace*{-.85em}(\num{62})  & \num{2}                                      & \num{2}  & \hspace*{-.85em}(\num{200})        & \num{2}          & \num{2}  & \hspace*{-.85em}(\num{200}) & \num{2}  \\
		3                                  & \num{2146689}                        & \num{2}  & \hspace*{-.85em}(\num{87})  & --                                           & \num{6}  & \hspace*{-.85em}(\num{600})        & --               & \num{6}  & \hspace*{-.85em}(\num{600}) & \num{6}  \\
		4                                  & \num{16974593}                       & \num{2}  & \hspace*{-.85em}(\num{137}) & --                                           & \num{4}  & \hspace*{-.85em}(\num{400})        & --               & \num{4}  & \hspace*{-.85em}(\num{400}) & --       \\
		5                                  & \num{135005697}                      & \num{2}  & \hspace*{-.85em}(\num{200}) & --                                           & \num{3}  & \hspace*{-.85em}(\num{300})        & --               & \num{3}  & \hspace*{-.85em}(\num{300}) & --       \\
		\bottomrule
	\end{tabular}
\end{table}


\subsection{Goal-oriented Adaptivity Driven by a Linear Functional}
\label{ex:regular, linear functional}
\footnotetext{\url{https://www.oeaw.ac.at/ricam/hpc}}
Next, we 
consider the same setting as in the previous Example~\ref{ex:regular, smooth}. 
However, we are now not interested in behaviour of the solution $u$ in 
the complete space-time cylinder $Q$, but only in the integral over the spatial domain $\Omega$
at final time $T$, i.e. we are interested in the goal functional 
\begin{equation*}
\label{eqn:GoalFunctional1}
J(u) = \int_{\Omega}\! u(\,.\,,T)\;\mathrm{d}\Omega. 
\end{equation*}
Since the exact solution $u$ is given by \eqref{eqn:Example5.1:ExactSolution},
we can compute the value of the goal functional
\[ 
	J(u) = \frac{4\,e}{\pi^2} \approx 1.10167812933171
\]
at the exact solution $u$.

\begin{figure}[!htb]
	\centering
	\includegraphics[width=.95\linewidth]{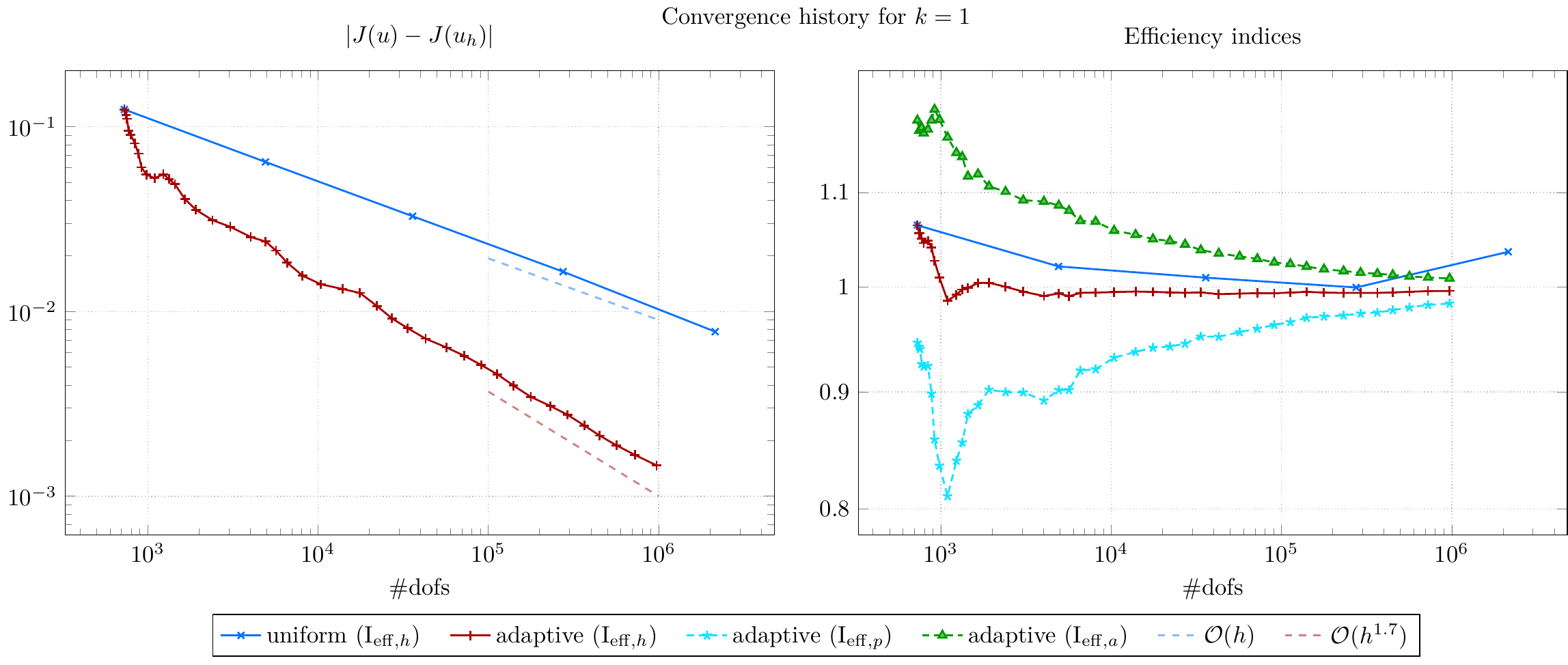}
	\caption{Example~\ref{ex:regular, linear functional} ($d=2$): Convergence history of the error in the functional 
	as well as efficiency plots, where we additionally included the efficiency of the primal and adjoint parts, respectively.}
	\label{fig2}
\end{figure}

\begin{figure}[!htb]
	\centering
	\includegraphics[width=.95\linewidth]{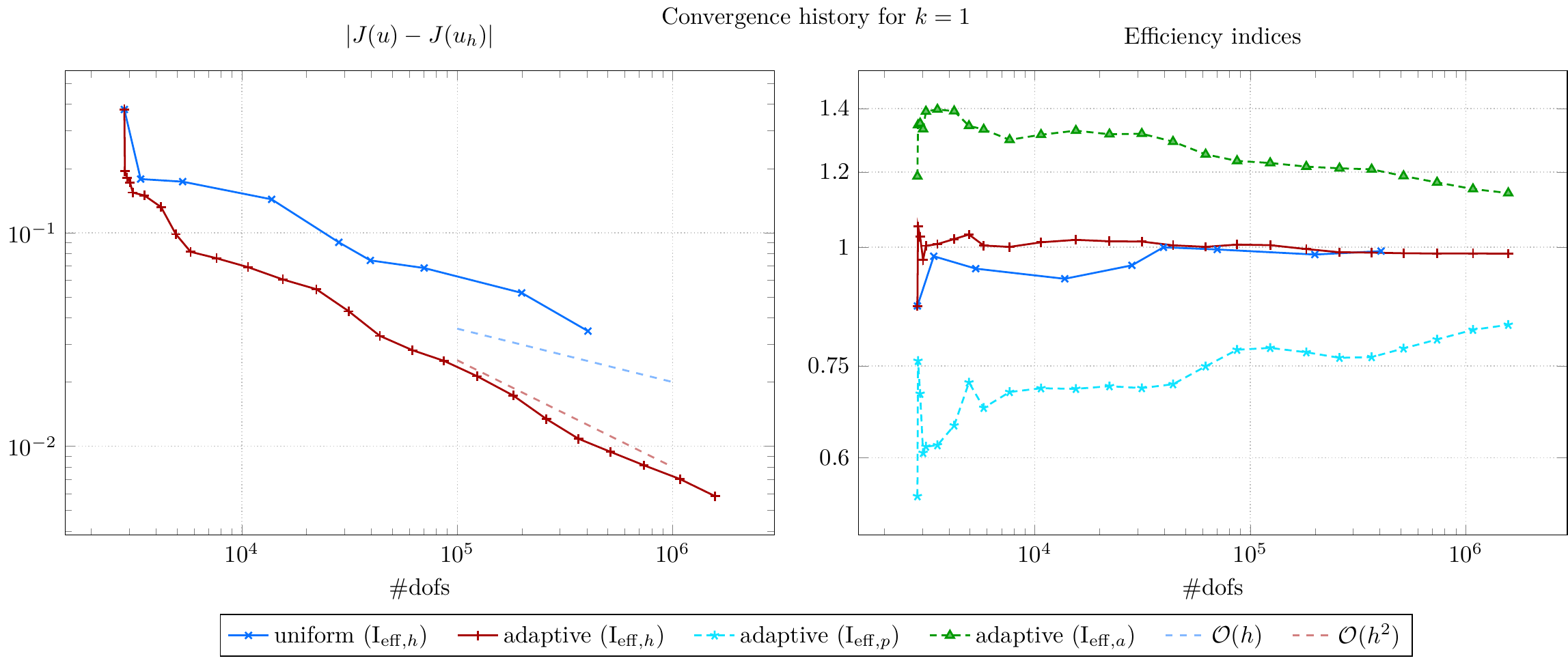}
	\caption{Example~\ref{ex:regular, linear functional} ($d=3$): Convergence history of the error in the functional 
	as well as efficiency plots, where we additionally included the efficiency of the primal and adjoint parts, respectively.}
	\label{ex2_4d}
\end{figure}

\begin{figure}[!htb]
	\centering%
	\includegraphics[width=.5\linewidth]{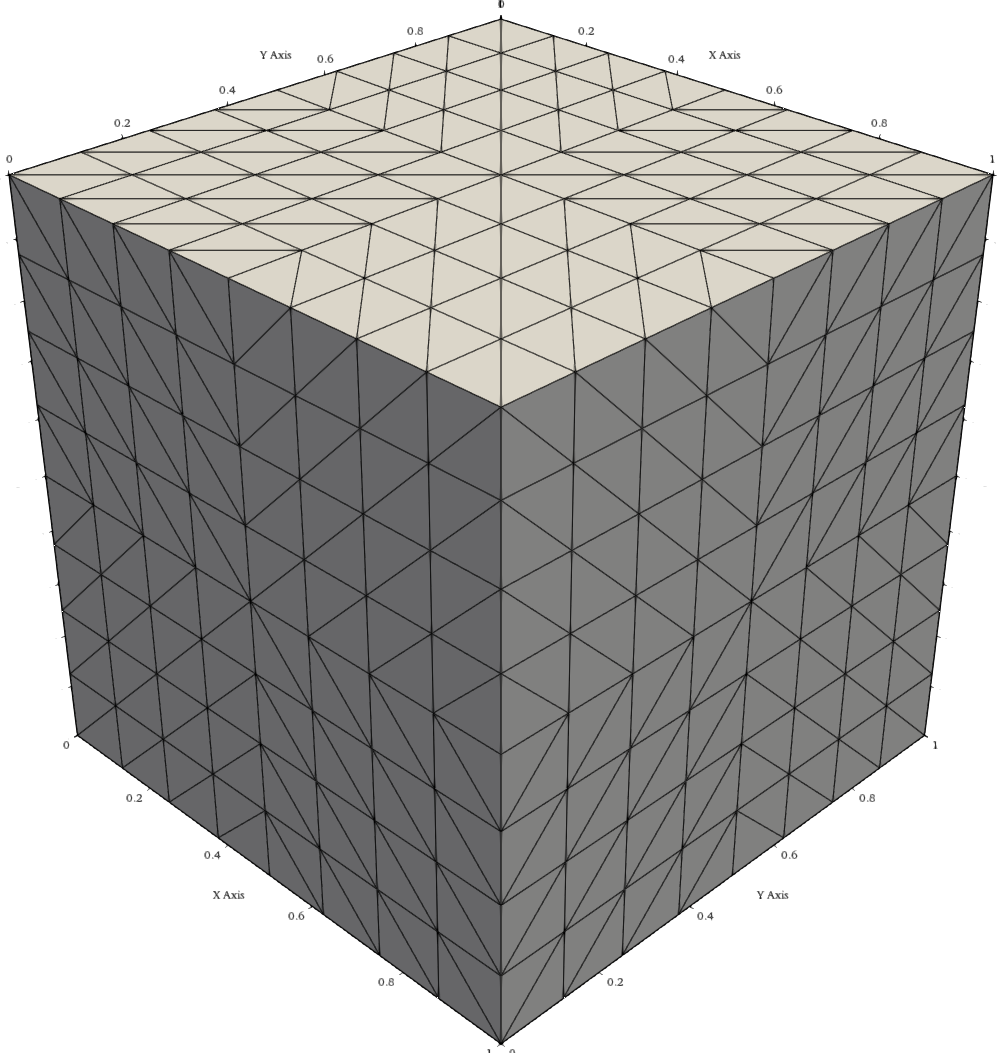}%
	\hfill%
	\includegraphics[width=.5\linewidth]{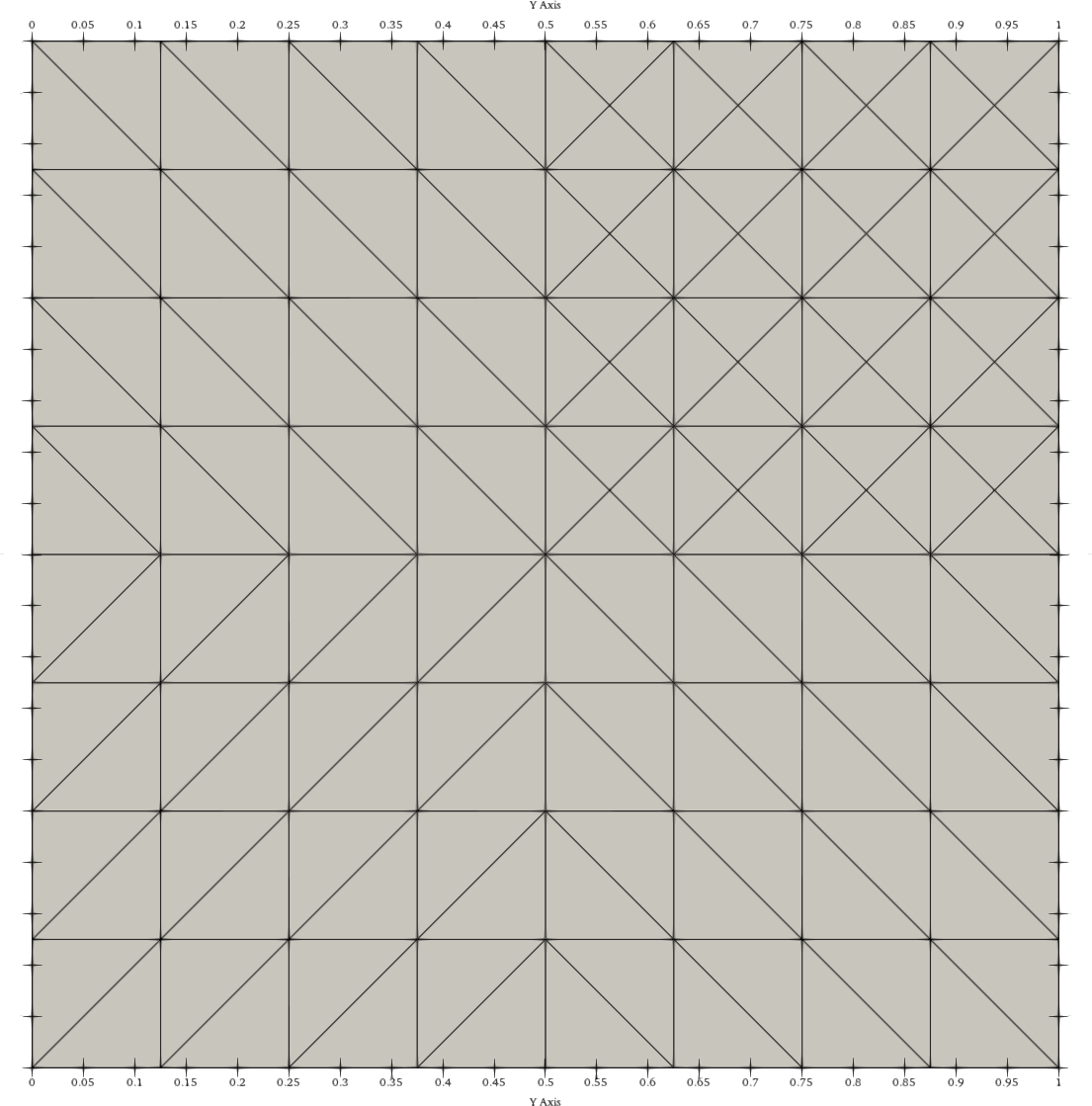}%
	\hfill%
	\includegraphics[width=.5\linewidth]{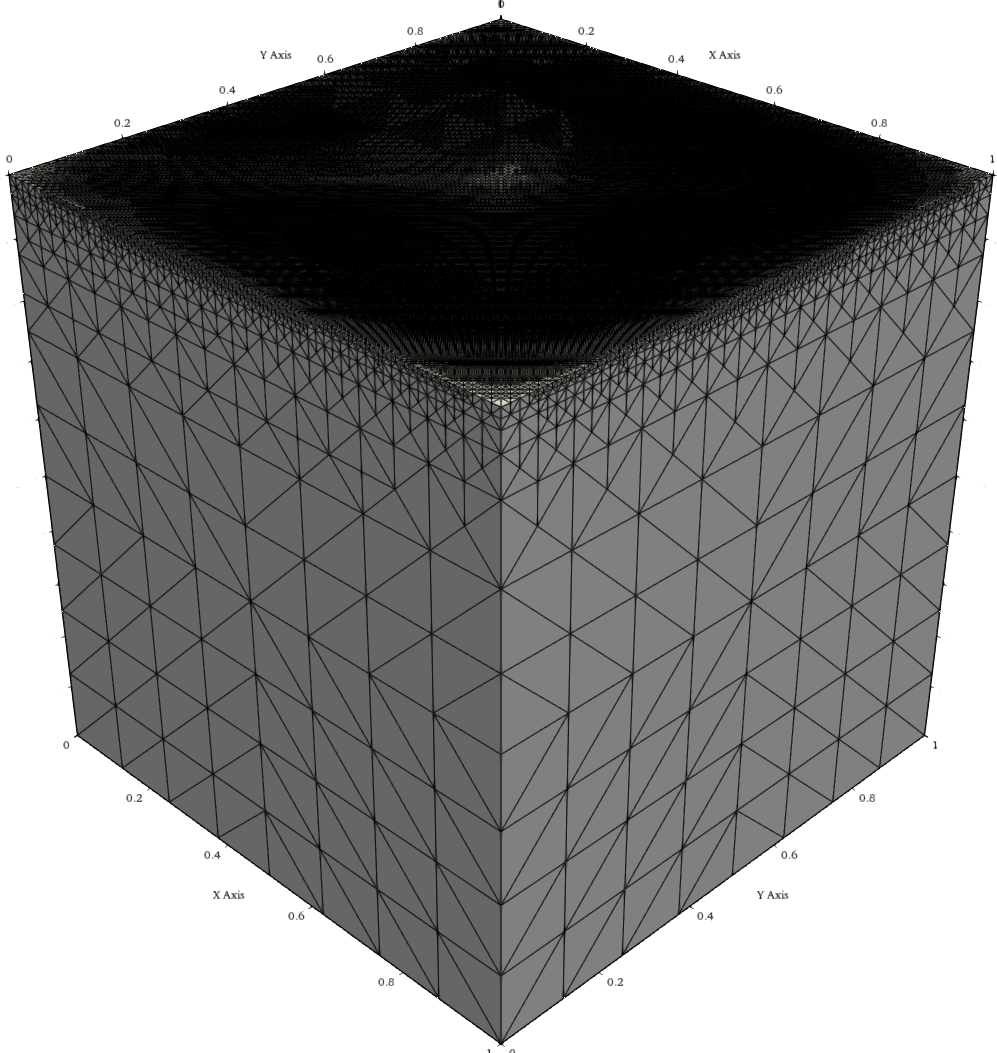}%
	\hfill%
	\includegraphics[width=.5\linewidth]{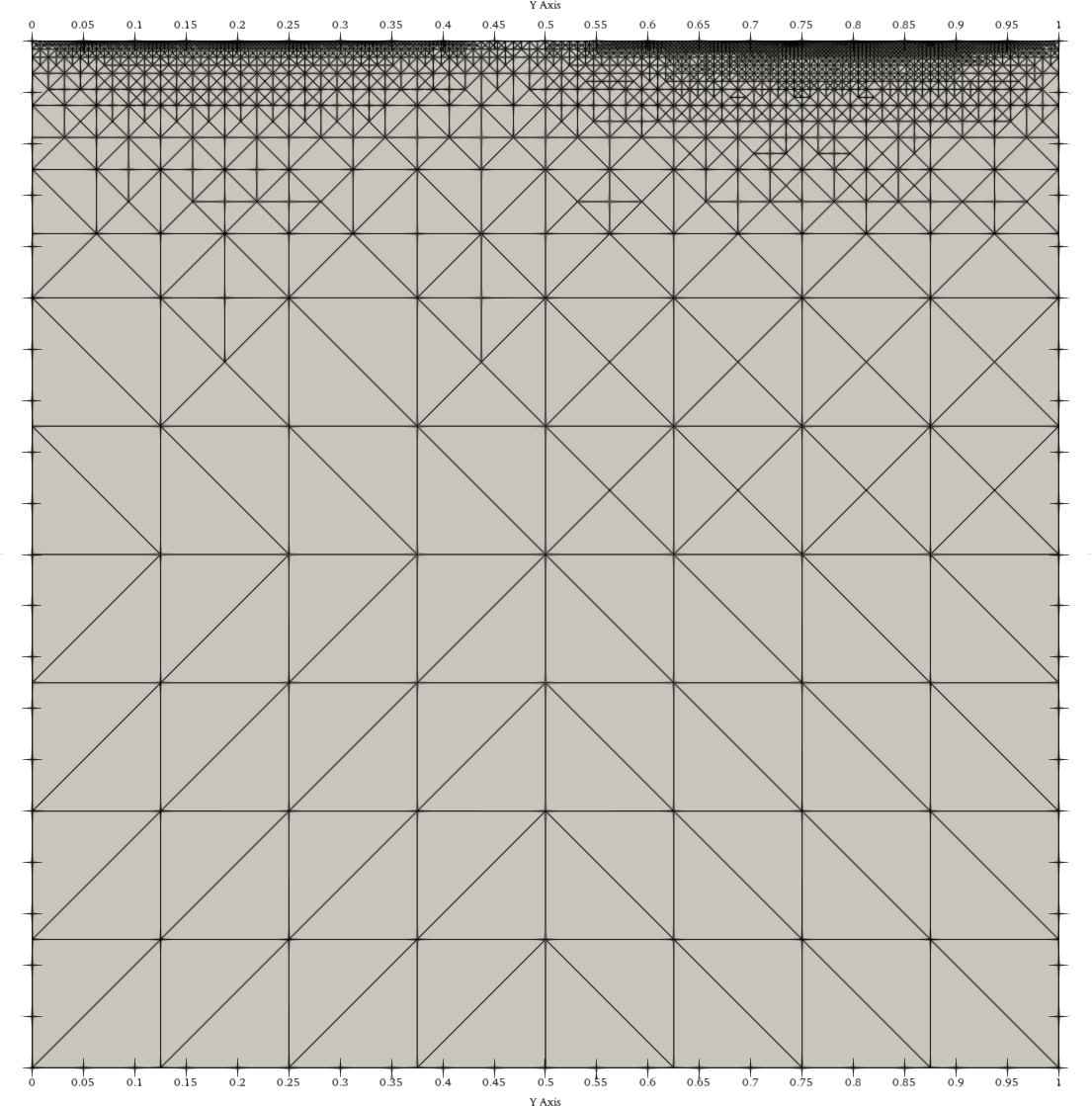}%
	\caption{Example~\ref{ex:regular, linear functional} ($d=2$): Initial space-time mesh (upper left); $(x_2,t)$-plane at $x_1 = 0.5$ of the initial space-time mesh (upper right); space-time mesh after 45 adaptive refinements (lower left); $(x_2,t)$-plane at $x_1 = 0.5$ after 45 adaptive refinements (lower right); using linear finite elements.}
	\label{fig3}
\end{figure}

Let us first consider the convergence history of the error in the linear functional $ J(\cdot)$
for $d=2$; see the left plot of Figure~\ref*{fig2}.
Here we observe that uniform refinements result in an error rate of $ O(h)$, where the mesh parameter $h$ is defined as $h = N_h^{-1/d}$, with $ N_h $ the total number of space-time dofs. 
On the other hand, using adaptive refinements driven by the goal-oriented error estimator, 
we obtain 
an improved
rate of $ O(h^{1.7}) $. Moreover, we also take a look at the efficiency index $\mathrm{I_{eff}}$ of the error estimator; cf. the right plot of Figure~\ref{fig2}. 
We observe that, after some initial oscillations, the efficiency index for the adaptive refinements remains close to $1$.
Figure~\ref{ex2_4d} shows that the same observations can be made in case $d=3$, 
where the  space-time cylinder $Q \subset \mathbb{R}^4$.
Next, let us consider the meshes produced by the adaptive finite element method. 
Since our quantity of interest is concentrated at the final time $T$, 
we expect heavy refinements towards the top of the space-time cylinder $Q$. 
In Figure~\ref{fig3}, we present the initial mesh as well as the mesh after a certain number of adaptive refinements
for the case $d=2$. 
Indeed, as we can observe in the lower row of Figure~\ref{fig3}, the mesh refinements are concentrated 
 
towards the top $\Sigma_T$ of the space-time cylinder $Q$.
This behaviour is also visible if we cut the space-time cylinder $Q$ along the $(x_2,t)$-plane at $x_1 = 0.5$; 
cf. the lower right plot of Figure~\ref{fig3}. 
The refinements seem to be entirely concentrated in the final quarter of the time interval $(0,T)$.  
This rather broad refinement is due to the mesh refinement algorithms that on the one hand isotropically refine each simplex, and on the other hand also refine the neighborhood of an element in order to prevent mesh degeneration. 

Let us once more consider the performance of the nested Newton method. 
As for the previous Example~\ref{ex:regular, smooth}, we investigate the influence of the regularization parameter $\varepsilon$, the inner solver, and the adaptive meshes on the solution process. 
In Table~\ref{tab:regular, linear functional:h-scaling}, we present the scaling 
with respect to the adaptive refinement levels $\ell$ for different regularization parameters $\varepsilon$.
At a first glance, the damped Newton's method seems rather robust wrt the parameters.

\begin{table}[!ht]
    \centering%
    \caption{Example~\ref{ex:regular, linear functional}: Scaling of the (damped) Newton solver with total number of inner solves in brackets, for $d+1=3$ and $k=1$.}%
   \label{tab:regular, linear functional:h-scaling}
    \begin{tabular}{clrlcrlcrlc}
        \toprule \multirow{2}{*}{$ \ell $} & \multirow{2}{*}{$ \#\mathrm{dofs} $} & \multicolumn{3}{c}{$ \varepsilon = 1$} & \multicolumn{3}{c}{$ \varepsilon = 10^{-5}$} & \multicolumn{3}{c}{$ \varepsilon = 10^{-10}$} \\
        & & \multicolumn{2}{c}{\textbf{GMRES}} & {\textbf{MUMPS}} & \multicolumn{2}{c}{\textbf{GMRES}} & {\textbf{MUMPS}} & \multicolumn{2}{c}{\textbf{GMRES}} & {\textbf{MUMPS}} \\ \midrule %
        \num{0} & \num{\sim729} & \num{9} & \hspace*{-0.85em}(\num{170}) & \num{9} & \num{11} & \hspace*{-0.85em}(\num{231}) & \num{11} & \num{11} & \hspace*{-0.85em}(\num{231}) & \num{11} \\
        \num{12} & \num{\sim4346} & \num{3} & \hspace*{-0.85em}(\num{112}) & \num{3} & \num{3} & \hspace*{-0.85em}(\num{171}) & \num{3} & \num{3} & \hspace*{-0.85em}(\num{171}) & \num{3} \\
        \num{20} & \num{\sim33886} & \num{2} & \hspace*{-0.85em}(\num{156}) & \num{2} & \num{2} & \hspace*{-0.85em}(\num{200}) & \num{2} & \num{2} & \hspace*{-0.85em}(\num{200}) & \num{2} \\
        \num{28} & \num{\sim264942} & \num{2} & \hspace*{-0.85em}(\num{200}) & \num{2} & \num{2} & \hspace*{-0.85em}(\num{200}) & \num{3} & \num{2} & \hspace*{-0.85em}(\num{200}) & \num{3} \\
        \num{34} & \num{\sim1112168} & \num{2} & \hspace*{-0.85em}(\num{200}) & -- & \num{2} & \hspace*{-0.85em}(\num{200}) & -- & \num{2} & \hspace*{-0.85em}(\num{200}) & -- \\
        \bottomrule
    \end{tabular}
\end{table}

	
\subsection{Goal-Oriented with a Non-linear Goal Functional}
\label{ex:nonlinear goal}

As our third example, we consider again the setting from the previous two examples. 
However, 
we are now interested in 
the following, non-linear volume goal functional
\[
	J(u) = \int_{Q_I}\! | \nabla_x u |^p\; dQ \approx 0.01937125060566419,
\]
where $ Q_I $ is a 
prescribed
region of interest. In our case, we choose $ Q_I $ to be an octahedron with edge length $0.5$, centered at $ (0.5,0.5,0.5) $; 
see 
also
the middle column of Figure~\ref{fig:nonlinear goal:meshes}.
We make this choice in order to exactly capture the region of interest with the finite element mesh. 
\begin{figure}[!htb]
	\centering
	\includegraphics[width=.95\linewidth]{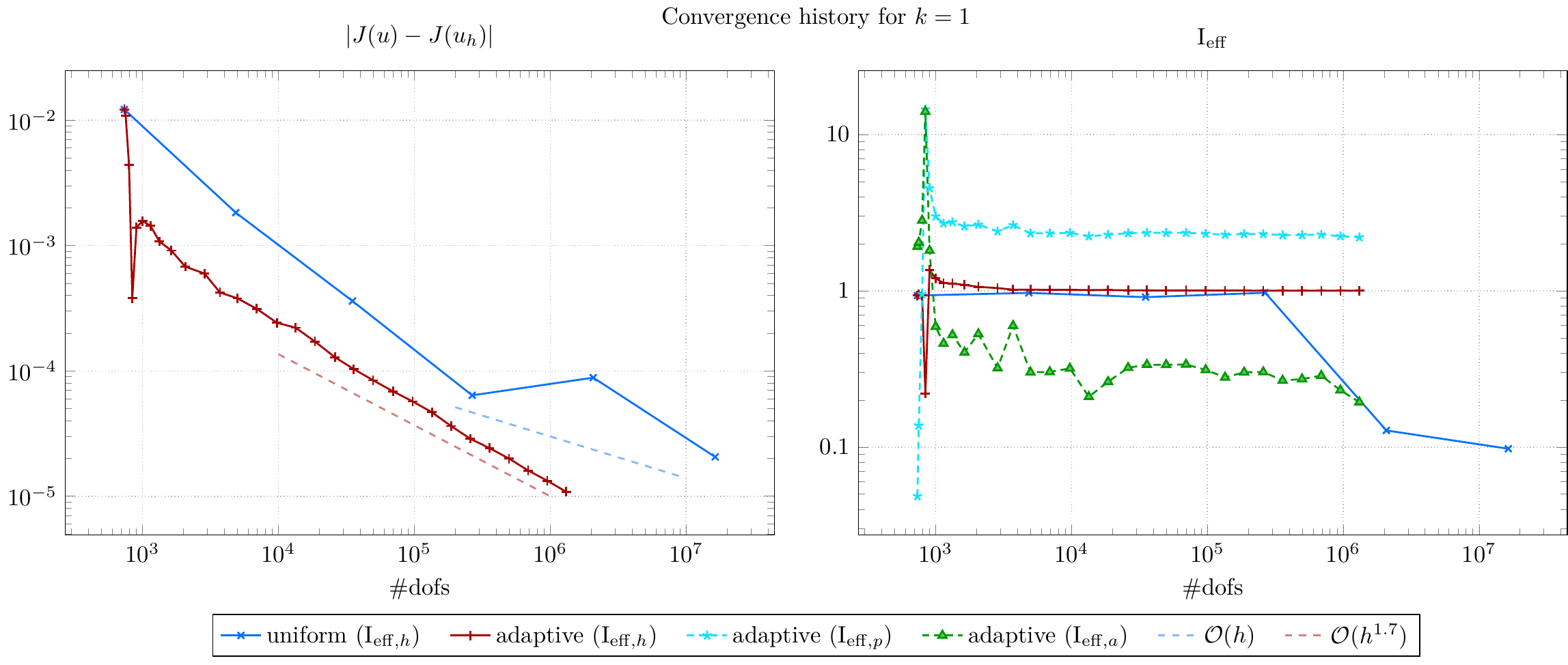}
	\caption{Example~\ref{ex:nonlinear goal} ($d=2$): Convergence history of the error in the functional 
	as well as efficiency plots, where we additionally included the efficiency of the primal and adjoint parts, respectively.}%
	\label{fig:nonlinear goal:convergence}
\end{figure}
Let us first 
consider the convergence history. In the left plot Fig.~\ref{fig:nonlinear goal:convergence}, we present the convergence history of the discretization error in the nonlinear goal functional $J(u)$. We observe that while uniform refinements reduce the overall error, adaptive refinements driven by the goal oriented error estimator lead to a considerable reduction in the number of dofs needed to attain a similar error. Moreover, we observe that the efficiency index of the adaptive refinements once more converges towards $ 1 $; cf. the right plot of Fig.~\ref{fig:nonlinear goal:convergence}.
\begin{figure}[!htb]
	\centering%
	\includegraphics[width=.3\linewidth]{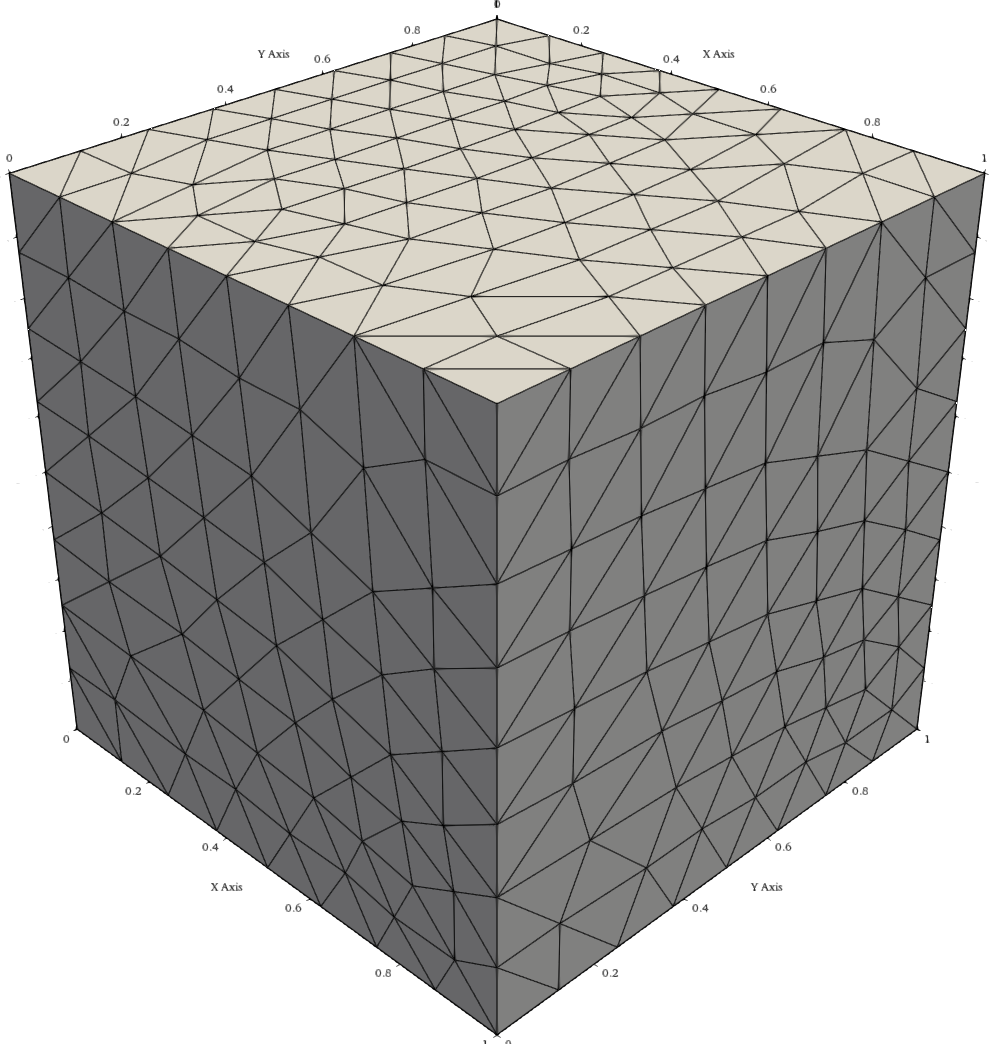}%
	\hfill%
	\includegraphics[width=.3\linewidth]{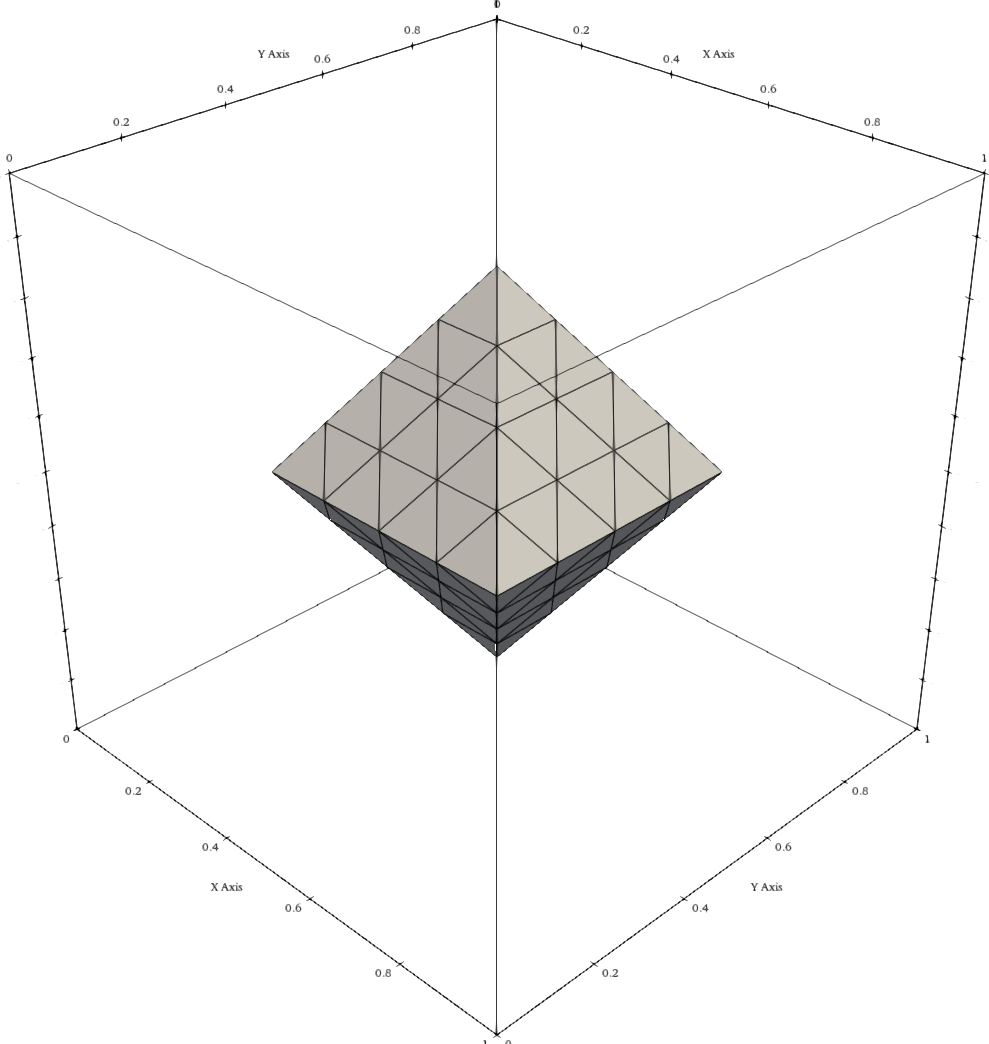}%
	\hfill%
	\includegraphics[width=.3\linewidth]{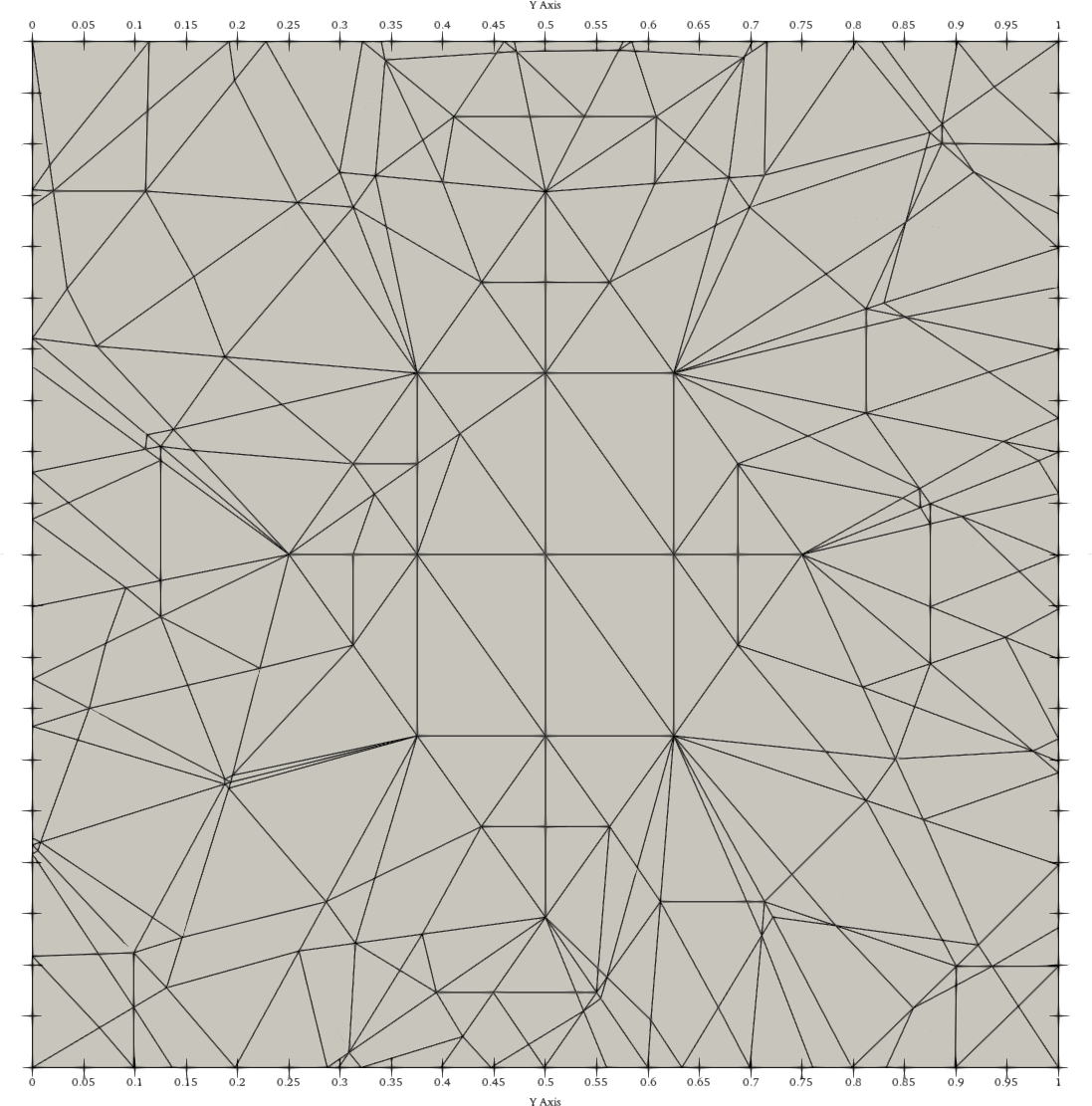}%
	\hfill%
	\includegraphics[width=.3\linewidth]{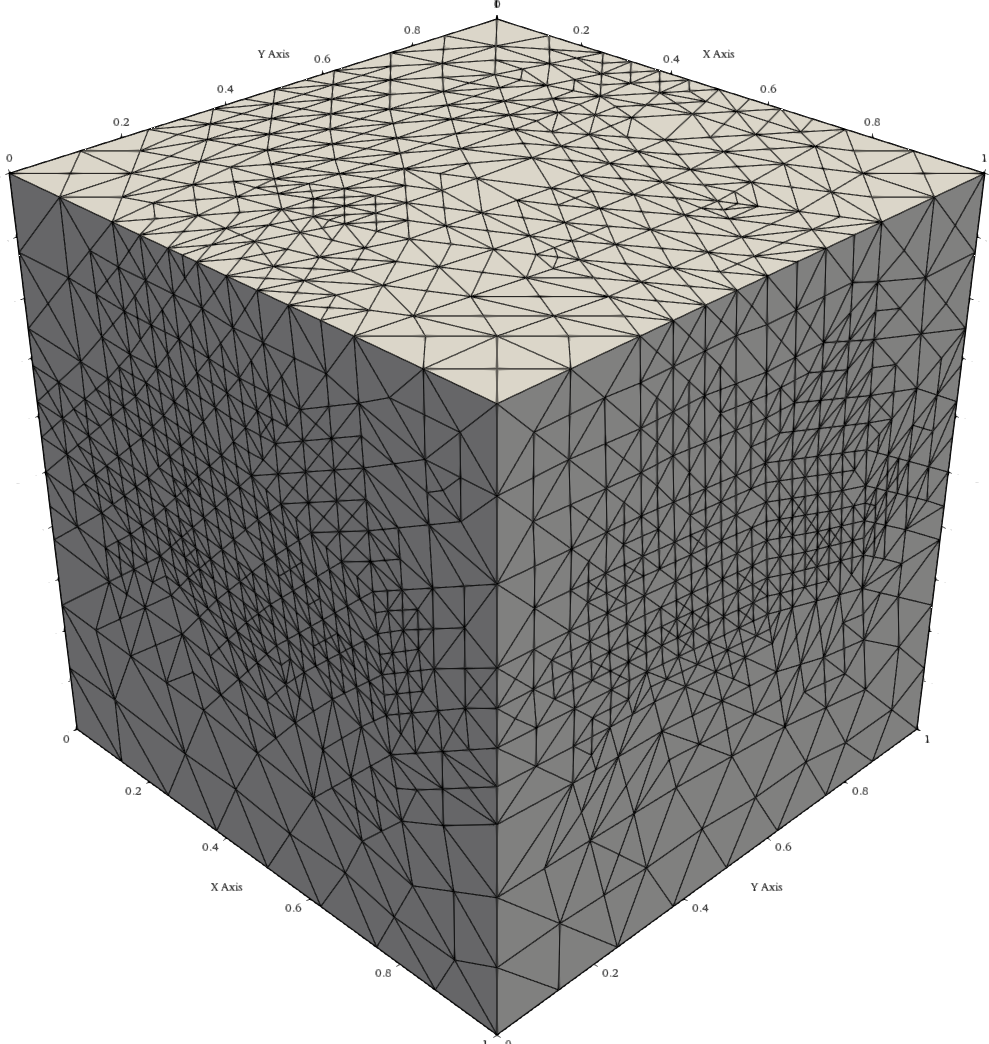}%
	\hfill%
	\includegraphics[width=.3\linewidth]{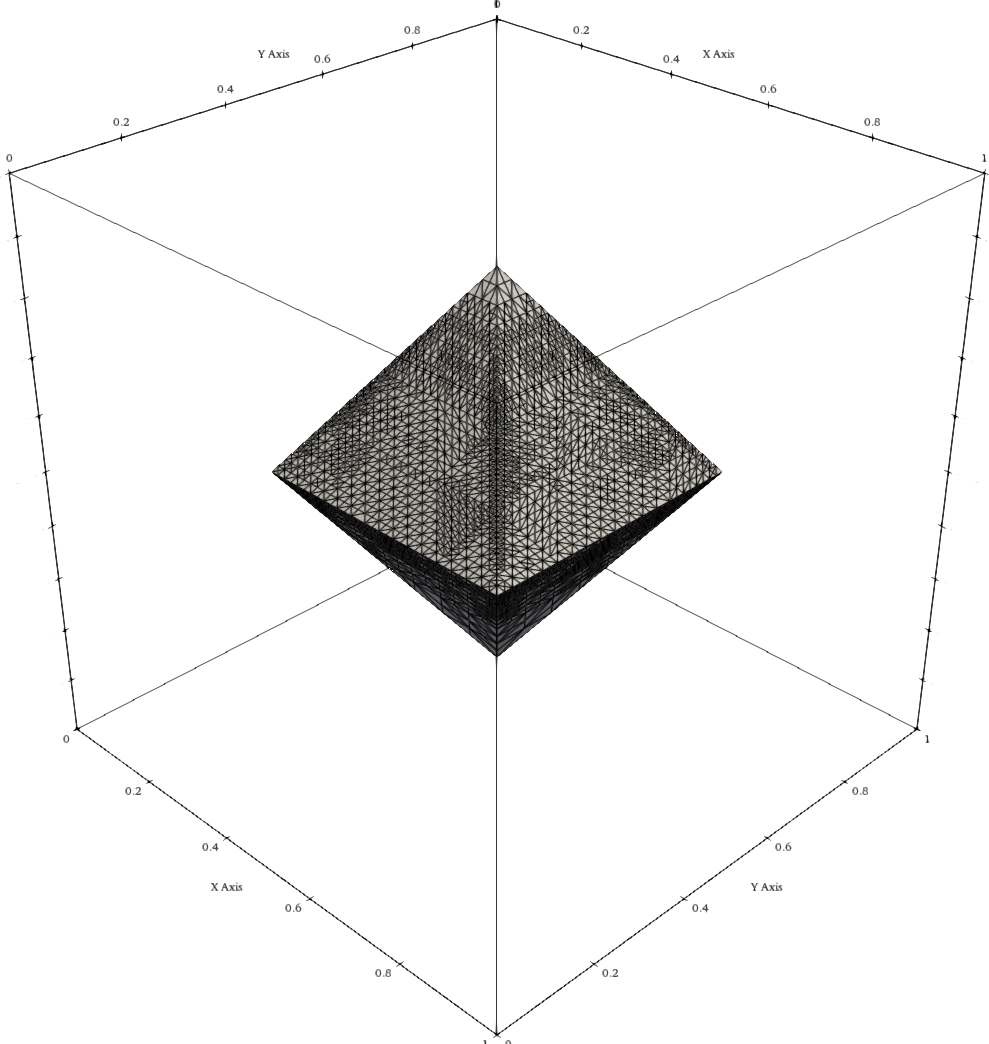}%
	\hfill%
	\includegraphics[width=.3\linewidth]{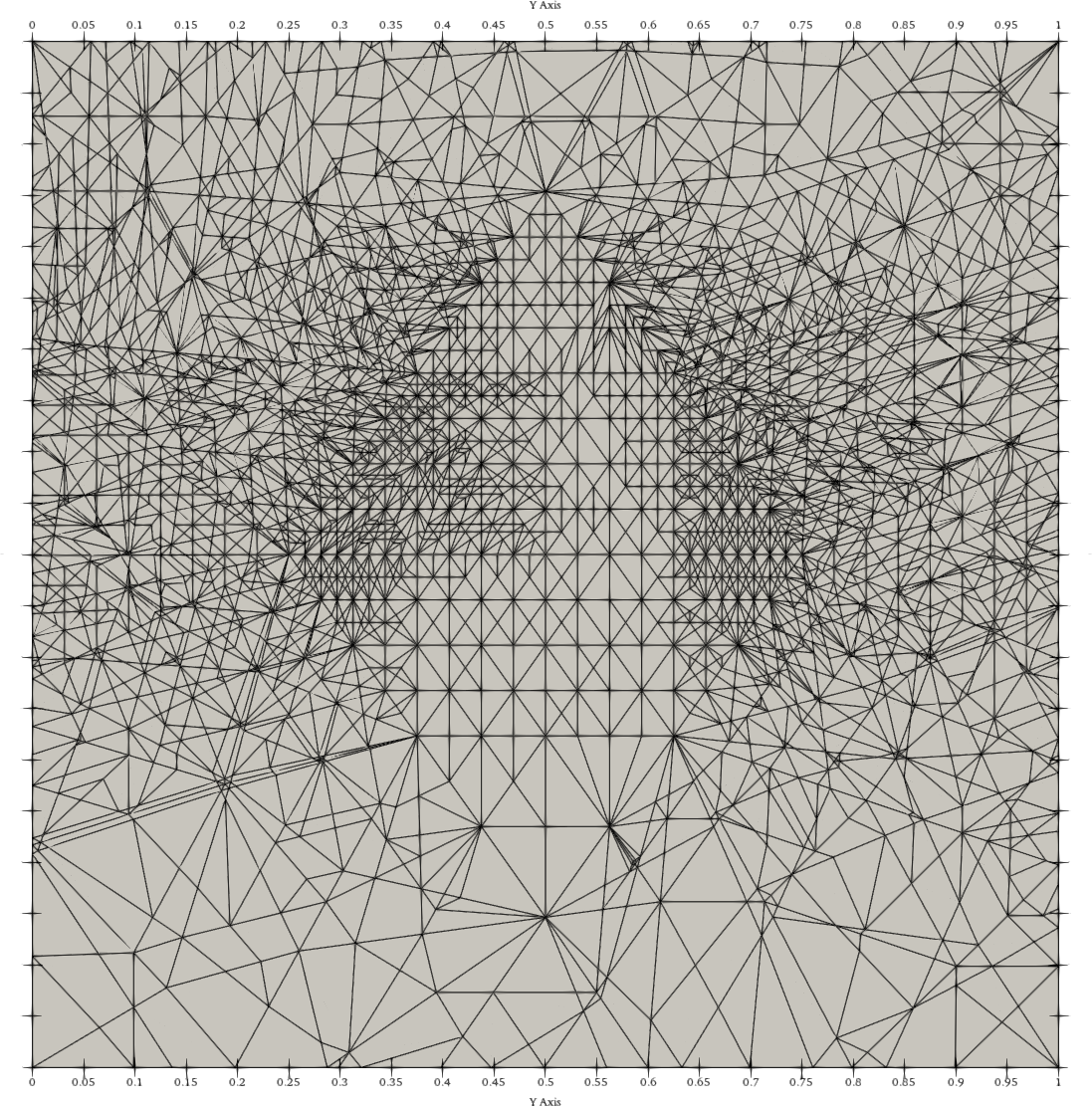}%
	\caption{Example~\ref{ex:nonlinear goal} ($d=2$): From left to right: full space-time mesh, surface mesh of $ Q_I$, and the $(x_2,t)$-plane at $x_1 = 0.5$; in its initial configuration (upper row), and after 20 adaptive refinements (lower row); using linear finite elements.}%
	\label{fig:nonlinear goal:meshes}
\end{figure}
The upper row of Figure~\ref{fig:nonlinear goal:meshes} presents the initial configuration of the finite element meshes, whereas the lower row shows the meshes after 20 adaptive refinements driven by the goal oriented estimator.
In the very left column, we can observe that the adaptive refinement indeed 
produces mostly unstructured space-time mesh.
Moreover, we see that the refinements are concentrated inside the octahedron $Q_I$, and any refinements outside are necessary in order to avoid any degeneration of the mesh elements. 

In Table~\ref{tab:nonlinear goal:h-scaling}, we again present the scaling wrt the level of refinement $\ell$ of the nested Newton method for different regularization parameters $ \varepsilon $. We again observe stable iteration counts between two and three Newton iterations as well as increasing iteration counts for the total number of inner solves for the Jacobian.

\begin{table}[!ht]
    \centering%
    \caption{Example~\ref{ex:nonlinear goal}: Scaling of the (damped) Newton solver with total number of inner solves in brackets, for $d+1=3$ and $k=1$.}
    \label{tab:nonlinear goal:h-scaling}
\begin{tabular}{clrlcrlcrlc}
    \toprule \multirow{2}{*}{$ \ell $} & \multirow{2}{*}{$ \#\mathrm{dofs} $} & \multicolumn{3}{c}{$ \varepsilon = 1$} & \multicolumn{3}{c}{$ \varepsilon = 10^{-5}$} & \multicolumn{3}{c}{$ \varepsilon = 10^{-10}$} \\
    & & \multicolumn{2}{c}{\textbf{GMRES}} & {\textbf{MUMPS}} & \multicolumn{2}{c}{\textbf{GMRES}} & {\textbf{MUMPS}} & \multicolumn{2}{c}{\textbf{GMRES}} & {\textbf{MUMPS}} \\ \midrule %
    \num{0} & \num{\sim735} & \num{10} & \hspace*{-0.85em}(\num{199}) & \num{10} & \num{11} & \hspace*{-0.85em}(\num{240}) & \num{11} & \num{11} & \hspace*{-0.85em}(\num{240}) & \num{11} \\
    \num{8} & \num{\sim4931} & \num{3} & \hspace*{-0.85em}(\num{81}) & \num{3} & \num{3} & \hspace*{-0.85em}(\num{104}) & \num{3} & \num{3} & \hspace*{-0.85em}(\num{104}) & \num{3} \\
    \num{15} & \num{\sim56031} & \num{2} & \hspace*{-0.85em}(\num{92}) & \num{2} & \num{2} & \hspace*{-0.85em}(\num{160}) & \num{2} & \num{2} & \hspace*{-0.85em}(\num{160}) & \num{2} \\
    \num{18} & \num{\sim155318} & \num{2} & \hspace*{-0.85em}(\num{114}) & \num{2} & \num{2} & \hspace*{-0.85em}(\num{200}) & \num{3} & \num{2} & \hspace*{-0.85em}(\num{200}) & \num{3} \\
    \num{24} & \num{\sim1133165} & \num{2} & \hspace*{-0.85em}(\num{200}) & -- & \num{3} & \hspace*{-0.85em}(\num{300}) & -- & \num{3} & \hspace*{-0.85em}(\num{300}) & -- \\
\bottomrule
\end{tabular}
\end{table}

\section{Conclusion and Outlook}
\label{Section:ConclustionOutlook}

We have proposed a new goal-oriented adaptive space-time finite element method 
for regularized parabolic $p$-Laplace initial-boundary value problems.
The space-time finite element discretization is based on the decomposition 
of the space-time cylinder into conforming simplicial elements like 
in the case of elliptic boundary value problems. 
The mesh refinement is driven by the DWR method and their localization 
is done
by means of the PU technique. So we can adaptively generate 
finite element approximations that are tailored to the quantity of interest (goal)
that is mathematically given by some possibly nonlinear functional.
Since we used a manufactured solution in our numerical experiments, 
we were able to compute the efficiency indices as quotient of the estimated 
error and the real error of the approximations to the functional. 
In all cases  
the efficiency indices were close to one. 
The adaptive process always saved a lot of unknowns in order to obtain some 
accuracy in the approximation of the functional in comparison with uniform refinement.
However, at each refinement level, we have to solve one non-linear system 
for the finite element solution and one linear linear system for the 
adjoint finite element solution. Furthermore, we need improved approximations
to primal and adjoint solutions that can be obtained by different methods 
as discussed in 
Subsection~\ref{Subsection:AnErrorIdentity}. 

A priori discretization error estimates as were proved in the case of 
linear parabolic initial-boundary value problem \cite{Steinbach:2015a},
convergence analysis of the adaptive process,
investigation of the convergence of the Newton solver,
improvement of the inner solver respectively preconditioner 
for the Jacobian system at each Newton iteration step,
and the adaption of the inner iteration to the convergence 
of the Newton iteration 
are future research topics. 
The latter topic as well as the simultaneous parallelization in space and time 
can lead to a considerably improvement 
of the performance of the algorithm.
This is also important 
for
more complex practical applications like
non-Newtonian flow problems described by power law models;
see e.g.\ \cite{Toulopoulos:2022b} and the references therein.

\section*{Acknowledgments}

This work has been supported by Deutsche Forschungsgemeinschaft (DFG, German Research Foundation)
under Germany's Excellence Strategy within the Cluster of
Excellence PhoenixD (EXC 2122), 
and 
by the  Austrian Science Fund (FWF) under the grant DK W1214-04.
Furthermore, 
the first author Bernhard Endtmayer
greatly acknowledges the support and funding of the 'Alexander von Humboldt Foundation' as well as their 'Humboldt Fellowship'.


		\bibliography{lit.bib}
		\bibliographystyle{abbrv}

\end{document}